\documentclass{article}%
\usepackage{bbm}
\usepackage{epsfig}
\usepackage{graphics,graphicx,amssymb,amsmath,verbatim}
\usepackage{amssymb}
\usepackage{mathrsfs}
\usepackage{amsfonts}
\usepackage{amsmath}
\usepackage{graphicx}
\usepackage{easybmat}
\usepackage{color}%
\providecommand{\U}[1]{\protect \rule{.1in}{.1in}}
\newtheorem{theorem}{Theorem}

\newtheorem{corollary}{Corollary}

\newtheorem{lemma}{Lemma}

\newtheorem{remark}{Remark}

\allowdisplaybreaks[4]

\begin{document}

\title{Analysis and Design of Complex-Valued Linear Systems}
\author{Bin Zhou\thanks{Center for Control Theory and Guidance Technology, Harbin
Institute of Technology, Harbin, 150001, China. Email: binzhoulee@163.com,
binzhou@hit.edu.cn.}}
\date{}
\maketitle

\begin{abstract}
This paper studies a class of complex-valued linear systems whose state
evolution dependents on both the state vector and its conjugate. The
complex-valued linear system comes from linear dynamical quantum control
theory and is also encountered when a normal linear system is controlled by
feedback containing both the state vector and its conjugate that can provide
more design freedom. By introducing the concept of bimatrix and its
properties, the considered system is transformed into an equivalent
real-representation system and a non-equivalent complex-lifting system,
which are normal linear systems. Some analysis and design problems including
solutions, controllability, observability, stability, eigenvalue assignment,
stabilization, linear quadratic regulation (LQR), and state observer design
are then investigated. Criterion, conditions, and algorithms are provided in
terms of the coefficient bimatrices of the original system. The developed
approaches are also utilized to investigate the so-called antilinear system
which is a special case of the considered complex-valued linear system. The
existing results on this system have been improved and some new results are
established.

\vspace{0.3cm}


\end{abstract}

\section{Introduction}

We study in this paper the following complex-valued linear system%
\begin{equation}
\left \{
\begin{array}{rl}
x^{+} & =A_{1}x+A_{2}^{\#}x^{\#}+B_{1}u+B_{2}^{\#}u^{\#}, \\
y & =C_{1}x+C_{2}^{\#}x^{\#}+D_{1}u+D_{2}^{\#}u^{\#},%
\end{array}%
\right.  \label{sys}
\end{equation}%
where $A_{i}\in \mathbf{C}^{n\times n},B_{i}\in \mathbf{C}^{n\times
m},C_{i}\in \mathbf{C}^{p\times n},D_{i}\in \mathbf{C}^{p\times m},i=1,2,$
are known coefficients, $x=x(t)$ is the state, $u=u(t)$ is the control, $%
y=y(t)$ is the output, $P^{\#}$ denotes the conjugate of the matrix/vector $%
P $, and $x^{+}(t)$ denotes $x\left( t+1\right) $ if $t\in \mathbf{Z}%
^{+}=\{0,1,2,\cdots \}$ (discrete-time systems) and denotes $\dot{x}(t)$ if $%
t\in \mathbf{R}^{+}=[0,\infty )$ (continuous-time systems). Throughout this
paper, the dependence of variables on $t$ will be suppressed unless
necessary. The initial condition is set to be $x\left( 0\right) =x_{0}\in
\mathbf{C}^{n}.$ Hereafter, the first equation of (\ref{sys}) is referred to
as the system equation while the second one is referred to as the output
equation.

Clearly, (\ref{sys}) becomes the normal linear system
\begin{equation}
\left \{
\begin{array}{rl}
x^{+} & =A_{1}x+B_{1}u, \\
y & =C_{1}x+D_{1}u,%
\end{array}%
\right.  \label{normal}
\end{equation}%
if $A_{2},B_{2},C_{2}$ and $D_{2}$ are null, and becomes the so-called
antilinear system
\begin{equation}
\left \{
\begin{array}{rl}
x^{+} & =A_{2}^{\#}x^{\#}+B_{2}^{\#}u^{\#}, \\
y & =C_{2}^{\#}x^{\#}+D_{2}^{\#}u^{\#},%
\end{array}%
\right.  \label{antilinear}
\end{equation}%
if $A_{1},B_{1},C_{1}$ and $D_{1}$ are absent. System (\ref{normal}) has
been very well studied in the literature (see, for example, %
\cite{kfa69book,rugh96book,xll17ijss,zyzy17ijss}) while system (\ref%
{antilinear}) without output equations were respectively studied in %
\cite{wdl13aucc} and \cite{wzls15iet}.

We study this class of complex-valued linear systems for several reasons.

\begin{itemize}
\item It is a natural extension of the well studied normal linear system (%
\ref{sys}) and the recently studied antilinear system (\ref{antilinear}).

\item In linear quantum systems theory, the Heisenberg evolution of the
annihilation--creation pairs and their conjugates are, in general, coupled
(see, for example, \cite{grn10pr}, \cite{zhang17auto}, \cite{zgpg18tac}, and pp. 44--49 in %
\cite{ny17book}). This leads to linear systems that are exactly in the form
of (\ref{sys}).

\item To stabilize the antilinear system (\ref{antilinear}) when $t\in
\mathbf{R}^{+},$ the so-called full state feedback
\begin{equation}
u=K_{1}x+K_{2}^{\#}x^{\#},  \label{eqfeedback}
\end{equation}%
is necessary (see the developments in Subsection \ref{sec4.1}), which
results in a closed-loop system in the form of (\ref{sys}).

\item Any real-valued normal linear system in the form of (\ref{normal}) can
be equivalently converted into the complex-valued linear system (\ref{sys})
whose dimension is only half.

\item Even for the normal linear system (\ref{normal}), the full state
feedback (\ref{eqfeedback}) is more powerful than the well-known normal
linear feedback%
\begin{equation}
u=K_{1}x,  \label{normalfeedback}
\end{equation}%
in the sense that (\ref{eqfeedback}) allows more design freedom (of course
we should disregard the physical implementability of the full state feedback
(\ref{eqfeedback})).
\end{itemize}

To illustrate the fourth reason, we consider a $2n$ dimensional normal
linear system in the form of (\ref{normal}) with real coefficients, namely,%
\begin{equation}
\left \{
\begin{array}{rl}
\xi ^{+} & =\left[
\begin{array}{cc}
A_{11} & A_{12} \\
A_{21} & A_{22}%
\end{array}%
\right] \xi +\left[
\begin{array}{cc}
B_{11} & B_{12} \\
B_{21} & B_{22}%
\end{array}%
\right] \upsilon , \\
\eta & =\left[
\begin{array}{cc}
C_{11} & C_{12} \\
C_{21} & C_{22}%
\end{array}%
\right] \xi +\left[
\begin{array}{cc}
D_{11} & D_{12} \\
D_{21} & D_{22}%
\end{array}%
\right] \upsilon ,%
\end{array}%
\right.  \label{eqsys6}
\end{equation}%
where $A_{ij}\in \mathbf{R}^{n\times n},B_{ij}\in \mathbf{R}^{n\times
m},C_{ij}\in \mathbf{R}^{p\times n},D_{ij}\in \mathbf{R}^{p\times m},i,j=1,2$
and $\xi \left( 0\right) =\xi _{0}\in \mathbf{R}^{2n}.$ Denote $\xi =[\xi
_{1}^{\mathrm{T}},\xi _{2}^{\mathrm{T}}]^{\mathrm{T}},\xi _{i}\in \mathbf{R}%
^{n},\upsilon =[\upsilon _{1}^{\mathrm{T}},\upsilon _{2}^{\mathrm{T}}]^{%
\mathrm{T}},\upsilon _{i}\in \mathbf{R}^{m},$ and $\eta =[\eta _{1}^{\mathrm{%
T}},\eta _{2}^{\mathrm{T}}]^{\mathrm{T}},\eta _{i}\in \mathbf{R}^{p},i=1,2.$
If we choose $x=\xi _{1}+\mathrm{j}\xi _{2},u=\upsilon _{1}+\mathrm{j}%
\upsilon _{2},$ and $y=\eta _{1}+\mathrm{j}\eta _{2},$ then system (\ref%
{eqsys6}) can be exactly written as the complex-valued linear system (\ref%
{sys}), where%
\begin{equation}
\left \{
\begin{array}{l}
A_{1}=\frac{1}{2}\left( A_{11}+A_{22}\right) +\frac{\mathrm{j}}{2}\left(
A_{21}-A_{12}\right) \in \mathbf{C}^{n\times n}, \\
A_{2}=\frac{1}{2}\left( A_{11}-A_{22}\right) -\frac{\mathrm{j}}{2}\left(
A_{21}+A_{12}\right) \in \mathbf{C}^{n\times n},%
\end{array}%
\right.  \label{eqadda1a2}
\end{equation}%
and $B_{i}\in \mathbf{C}^{n\times m},C_{i}\in \mathbf{C}^{p\times
n},D_{i}\in \mathbf{C}^{p\times m},i=1,2,$ are defined in a similar way. It
follows that the dimension of the corresponding complex-valued linear system
(\ref{sys}) is only half of that for system (\ref{eqsys6}). This is helpful,
at least, in simulation since the number of integrators has been halved.

\begin{remark}
We consider two special cases of system (\ref{eqsys6}). The first case is
the so-called symmetric linear system \cite{bodson14cdc,db16tie}%
\begin{equation}
\left \{
\begin{array}{rl}
\xi ^{+} & =\left[
\begin{array}{cc}
A_{11} & -A_{21} \\
A_{21} & A_{11}%
\end{array}%
\right] \xi +\left[
\begin{array}{cc}
B_{11} & -B_{21} \\
B_{21} & B_{11}%
\end{array}%
\right] \upsilon , \\
\eta & =\left[
\begin{array}{cc}
C_{11} & -C_{21} \\
C_{21} & C_{11}%
\end{array}%
\right] \xi +\left[
\begin{array}{cc}
D_{11} & -D_{21} \\
D_{21} & D_{11}%
\end{array}%
\right] \upsilon ,%
\end{array}%
\right.  \label{symmetricsys}
\end{equation}%
where $A_{ij}\in \mathbf{R}^{n\times n},B_{ij}\in \mathbf{R}^{n\times
m},C_{ij}\in \mathbf{R}^{p\times n},D_{ij}\in \mathbf{R}^{p\times
m},i,j=1,2. $ Such a system can be used to model many practical systems that
are symmetric in a sense, for example, the doubly-fed induction machine (see %
\cite{bodson14cdc,db16tie} and the references therein). Then it follows
from (\ref{eqadda1a2}) that the corresponding complex-valued linear system (%
\ref{sys}) is exactly in the form of (\ref{normal}) where $x=\xi _{1}+%
\mathrm{j}\xi _{2},y=\eta _{1}+\mathrm{j}\eta _{2},u=\upsilon _{1}+\mathrm{j}%
\upsilon _{2}$ and%
\begin{eqnarray*}
A_{1} &=&A_{11}+\mathrm{j}A_{21},B_{1}=B_{11}+\mathrm{j}B_{21},C_{1}=C_{11}+%
\mathrm{j}C_{21}, \\
D_{1} &=&D_{11}+\mathrm{j}D_{21},A_{2}=0,B_{2}=0,C_{2}=0,D_{2}=0.
\end{eqnarray*}%
The second case is the following so-called antisymmetric linear system%
\begin{equation}
\left \{
\begin{array}{rl}
\xi ^{+} & =\left[
\begin{array}{cc}
A_{11} & A_{21} \\
A_{21} & -A_{11}%
\end{array}%
\right] \xi +\left[
\begin{array}{cc}
B_{11} & B_{21} \\
B_{21} & -B_{11}%
\end{array}%
\right] \upsilon , \\
\eta & =\left[
\begin{array}{cc}
C_{11} & C_{21} \\
C_{21} & -C_{11}%
\end{array}%
\right] \xi +\left[
\begin{array}{cc}
D_{11} & D_{21} \\
D_{21} & -D_{11}%
\end{array}%
\right] \upsilon ,%
\end{array}%
\right.  \label{antisymmetricsys}
\end{equation}%
where $A_{ij}\in \mathbf{R}^{n\times n},B_{ij}\in \mathbf{R}^{n\times
m},C_{ij}\in \mathbf{R}^{p\times n},D_{ij}\in \mathbf{R}^{p\times
m},i,j=1,2. $ Such a system has a dual structure of (\ref{symmetricsys}),
and may have potential applications in engineering (thought were not found
so far). Then it follows from (\ref{eqadda1a2}) that the corresponding
complex-valued linear system (\ref{sys}) is exactly in the form of (\ref%
{antilinear}) where $x=\xi _{1}+\mathrm{j}\xi _{2},y=\eta _{1}+\mathrm{j}%
\eta _{2},u=\upsilon _{1}+\mathrm{j}\upsilon _{2}$ and%
\begin{eqnarray*}
A_{1} &=&0,B_{1}=0,C_{1}=0,D_{1}=0,A_{2}=A_{11}-\mathrm{j}A_{21}, \\
B_{1} &=&B_{11}-\mathrm{j}B_{21},C_{1}=C_{11}-\mathrm{j}C_{21},D_{1}=D_{11}-%
\mathrm{j}D_{21}.
\end{eqnarray*}
\end{remark}

To illustrate the final reason, we consider a normal linear system in the
form of (\ref{normal}) with%
\begin{equation}
A_{1}=\left[
\begin{array}{cc}
0 & 1 \\
-\alpha _{0} & -\alpha _{1}%
\end{array}%
\right] ,\;B_{1}=\left[
\begin{array}{c}
0 \\
1%
\end{array}%
\right] ,  \label{eqa1b1}
\end{equation}%
which is controllable, where $\alpha _{i}\in \mathbf{R},i=0,1$. Then, for
any desired characteristic polynomial $\gamma \left( s\right) =s^{2}+\gamma
_{1}s+\gamma _{0},\gamma _{i}\in \mathbf{R},i=0,1$, there is a unique $%
K_{1}\in \mathbf{R}^{1\times 2}$ such that the closed-loop system consisting
of (\ref{normal}) and (\ref{normalfeedback}) has the characteristic
polynomial $\gamma \left( s\right) $ \cite{kfa69book}. Such a unique $K_{1}$
is given by \cite{kfa69book}%
\begin{equation}
K_{1}=\left[
\begin{array}{cc}
\alpha _{0}-\gamma _{0} & \alpha _{1}-\gamma _{1}%
\end{array}%
\right] .  \label{eqk1exam}
\end{equation}%
Now if we allow the full state feedback (\ref{eqfeedback}), where $%
K_{1}=[\kappa _{11},\kappa _{12}]=[k_{111}+\mathrm{j}k_{112},k_{121}+\mathrm{%
j}k_{122}]$ and $K_{2}=[\kappa _{21},\kappa _{22}]=[k_{211}+\mathrm{j}%
k_{212},k_{221}+\mathrm{j}k_{222}]$ with $k_{ijl}\in \mathbf{R},$ the
closed-loop system is%
\begin{equation}
x^{+}=\left[
\begin{array}{cc}
0 & 1 \\
\kappa _{11}-\alpha _{0} & \kappa _{12}-\alpha _{1}%
\end{array}%
\right] x+\left[
\begin{array}{cc}
0 & 0 \\
\kappa _{21}^{\#} & \kappa _{22}^{\#}%
\end{array}%
\right] x^{\#},  \label{closed2}
\end{equation}%
which, by separating real and imaginary parts, is equivalent to%
\begin{equation}
\vec{x}^{+}=\left( A+BK\right) \vec{x},  \label{sys3}
\end{equation}%
where $\vec{x}=[\mathrm{Re}(x^{\mathrm{T}}),\mathrm{Im}(x^{\mathrm{T}})]^{%
\mathrm{T}},A=\mathrm{diag}\{A_{1},A_{1}\},B=\mathrm{diag}\{B_{1},B_{1}\},$
and%
\begin{equation}
K=\left[
\begin{array}{cccc}
k_{111}+k_{211} & k_{121}+k_{221} & -\left( k_{112}+k_{212}\right) & -\left(
k_{122}+k_{222}\right) \\
k_{112}-k_{212} & k_{122}-k_{222} & k_{111}-k_{211} & k_{121}-k_{221}%
\end{array}%
\right] .  \label{eqf}
\end{equation}%
As system (\ref{sys3}) is of order 4, we expand the associated desired
characteristic polynomial as $\gamma ^{2}\left( s\right) $ (notice that $%
\gamma \left( s\right) $ and $\gamma ^{2}\left( s\right) $ have the same
zeroes with only different multiplicities. Further explanation on such an
expansion will be given in Subsection \ref{sec2.2}). As $\left( A,B\right) $
is controllable and $B$ has two columns, for any $\gamma ^{2}\left( s\right)
,$ there exists non-unique matrix $K\in \mathbf{R}^{2\times 4}$ such that %
\cite{kfa69book}
\begin{equation}
\left \vert sI_{4}-\left( A+BK\right) \right \vert =\gamma ^{2}\left(
s\right) .  \label{eq54}
\end{equation}%
Denote any such a matrix $K$ by $[k_{ij}].$ Then it follows from (\ref{eqf})
that $[k_{11},k_{12},k_{13},k_{14},k_{21},k_{22},k_{23},k_{24}]^{\mathrm{T}%
}=T[k_{111},k_{112},k_{121},k_{122},k_{211},k_{212},k_{221},k_{222}]^{%
\mathrm{T}},$ where%
\begin{equation*}
T=\left[
\begin{array}{cccccccc}
1 & 0 & 0 & 0 & 1 & 0 & 0 & 0 \\
0 & 0 & 1 & 0 & 0 & 0 & 1 & 0 \\
0 & -1 & 0 & 0 & 0 & -1 & 0 & 0 \\
0 & 0 & 0 & -1 & 0 & 0 & 0 & -1 \\
0 & 1 & 0 & 0 & 0 & -1 & 0 & 0 \\
0 & 0 & 0 & 1 & 0 & 0 & 0 & -1 \\
1 & 0 & 0 & 0 & -1 & 0 & 0 & 0 \\
0 & 0 & 1 & 0 & 0 & 0 & -1 & 0%
\end{array}%
\right] .
\end{equation*}%
As $T$ is nonsingular, any $K$ determines uniquely a pair of $\left(
K_{1},K_{2}\right) ,$ which is not unique as $K$ is not unique. For
illustration purpose, we provide two different pairs of $\left(
K_{1},K_{2}\right) $ satisfying (\ref{eqf}) and (\ref{eq54}) as:%
\begin{align*}
\left[
\begin{array}{c}
K_{1} \\
K_{2}%
\end{array}%
\right] & =\left[
\begin{array}{cc}
f_{1}-\frac{\mathrm{j}}{2}\gamma _{0}^{2}-\frac{\mathrm{j}}{2} & \alpha
_{1}-\gamma _{1}-\mathrm{j}\gamma _{0}\gamma _{1} \\
f_{2}+\frac{\mathrm{j}}{2}\gamma _{0}^{2}-\frac{\mathrm{j}}{2} & \gamma _{1}+%
\mathrm{j}\gamma _{0}\gamma _{1}%
\end{array}%
\right] , \\
\left[
\begin{array}{c}
K_{1} \\
K_{2}%
\end{array}%
\right] & =\left[
\begin{array}{cc}
f_{1}+\frac{\mathrm{j}}{2}\gamma _{0}^{2}+\frac{\mathrm{j}}{2} & \alpha
_{1}-\gamma _{1}+\mathrm{j}\gamma _{0}\gamma _{1} \\
-f_{2}+\frac{\mathrm{j}}{2}\gamma _{0}^{2}-\frac{\mathrm{j}}{2} & -\gamma
_{1}+\mathrm{j}\gamma _{0}\gamma _{1}%
\end{array}%
\right] ,
\end{align*}%
in which $f_{1}=\alpha _{0}-\gamma _{0}-\frac{1}{2}\gamma _{1}^{2}$ and $%
f_{2}=\gamma _{0}+\frac{1}{2}\gamma _{1}^{2}.$ This is different from (\ref%
{eqk1exam}) which is unique.

For these reasons, in this paper we make a comprehensive study on the
complex-valued linear system (\ref{sys}). The considered problems include
solutions (state response), controllability, observability, stability,
eigenvalue assignment, stabilization, linear quadratic regulation (LQR), and
state observer design. Our study is based on two alternative descriptions of
system (\ref{sys}), which were motivated by literature for linear quantum
systems theory (for example, \cite{grn10pr,ny17book} and \cite{zhang17auto}%
), and the tool \textquotedblleft bimatrix\textquotedblright \ which has
many interesting and important properties to be studied in this paper. For
system analysis, we provide different equivalent criterion in terms of the
original system parameters, and for system design we provide necessary and
sufficient conditions to guarantee the existence of a solution, and give
appropriate algorithms whenever necessary. Most of these conditions are
expressed by the solvability of certain coupled equations involving the
bimatrices of the system, such as Lyapunov bimatrix equations and bimatrix
algebraic Riccati equations (AREs), whose coefficients are just the
coefficients for the original complex-valued linear system. We are
particularly interested in the analysis and design of the antilinear system (%
\ref{antilinear}). We show that all the coupled matrix equations encountered
in system analysis and design can be decoupled. By our approach we can not
only improve the existing results on this system, but also derive some new
results that were not available in the literature. For example, in the LQR
problem we are able to weaken the controllability condition in %
\cite{wqls16jfi} as stabilizability, and deal with both continuous and
discrete time antilinear systems. An abridged version of this paper was
presented on the 2017 Chinese Automation Congress (CAC) \cite{zhou17caa}.

The remainder of this paper is organized as follows. In Section \ref{sec2},
we introduce the notation of bimatrix as a fundamental tool for studying the
complex-valued linear system. System analysis including solutions,
controllability, observability and stability are carried out in Section \ref%
{sec3}, while system design including stabilization, eigenvalue assignment,
LQR, and state observer design are performed in Section \ref{sec4}. The
paper is finally concluded in Section \ref{sec5}.

\textbf{Notation}: We will use standard notation in this paper. For a matrix
$A\in \mathbf{C}^{n\times m},$we use $A^{\#},$ $A^{\mathrm{T}},$ $A^{\mathrm{%
H}},$ $\mathrm{rank}\left( A\right) ,$ $\left \vert A\right \vert ,$ $%
\left
\Vert A\right \Vert ,$ $\lambda \left( A\right) ,$ $\rho \left(
A\right) ,$ $\mathrm{Re}\left( A\right) $ and $\mathrm{Im}\left( A\right) $
to denote respectively its conjugate, transpose, conjugate transpose, rank,
determinant (when $n=m$), norm, eigenvalue set (when $n=m$), spectral radius
(when $n=m$), real part and imaginary part. Thus $A^{-\#}$ denotes $%
(A^{\#})^{-1}$ or $(A^{-1})^{\#}.$ We use $0_{n\times m}$ to denote an $%
n\times m$ zero matrix. For two integers $p,q$ with $p\leq q,$ we denote $%
\mathbf{I}\left[ p,q\right] =\{p,p+1,\cdots,q\}.$ Denote $\mathbf{Z}%
^{+}=\{0,1,2,\cdots \}$, $\mathbf{R}^{+}=[0,\infty),\mathbf{R}=\mathbf{R}%
^{+}\cup \{-\mathbf{R}^{+}\},\mathbf{Z}=\mathbf{Z}^{+}\cup \{-\mathbf{Z}%
^{+}\},$ $\otimes$ the Kronecker product, and $\mathrm{j}$ the unitary
imaginary number. For a series of matrices $A_{i},i\in \mathbf{I}\left[ 1,l%
\right] ,$ $\mathrm{diag}\{A_{1},A_{2},\cdots,A_{l}\}$ refers to a diagonal
matrix whose diagonal elements are $A_{i},i\in \mathbf{I}\left[ 1,l\right] .$
Other non-standard symbols will be defined when appear firstly.

\section{\label{sec2}Preliminaries}

\subsection{The bimatrix and its properties}

For a matrix pair $\left( A_{1},A_{2}\right) \in \left( \mathbf{C}^{n\times
m},\mathbf{C}^{n\times m}\right) $, denote the bimatrix $\left \{
A_{1},A_{2}\right \} $ that maps any $x\in \mathbf{C}^{m}\ $to $\mathbf{C}%
^{n} $ according to
\begin{equation}
y=\left \{ A_{1},A_{2}\right \} x\triangleq A_{1}x+A_{2}^{\#}x^{\#}.
\label{opeator}
\end{equation}%
If $x$ and $y$ are operators on Hilbert spaces, this is referred to as the coupling
operator in the quantum control literature (see \cite{grn10pr,ny17book} and the references therein). Define $%
(\left \{ A_{1},A_{2}\right \} +\left \{ B_{1},B_{2}\right \} )x=\left \{
A_{1},A_{2}\right \} x+\left \{ B_{1},B_{2}\right \} x$ and $\left \{
A_{1},A_{2}\right \} \left \{ B_{1},B_{2}\right \} x=\left \{
A_{1},A_{2}\right \} (\left \{ B_{1},B_{2}\right \} x).$ We denote $\left \{
A_{1},A_{2}\right \} =\left \{ B_{1},B_{2}\right \} \in \{ \mathbf{C}%
^{n\times m},\mathbf{C}^{n\times m}\}$ if $\left \{ A_{1},A_{2}\right \}
x=\left \{ B_{1},B_{2}\right \} x,\forall x\in \mathbf{C}^{m}.$ Then the
bimatrix has some properties listed below (the matrices involved are assumed
to have suitable dimensions):

\begin{enumerate}
\item $\{A_{1},A_{2}\} \left( x+y\right) =\left \{ A_{1},A_{2}\right \}
x+\left \{ A_{1},A_{2}\right \} y,\forall x,y\in \mathbf{C}^{m}$.

\item $\left \{ A_{1},A_{2}\right \} ax=a\left \{ A_{1},A_{2}\right \}
x,\forall x\in \mathbf{C}^{m},\forall a\in \mathbf{R}.$

\item $\left \{ A_{1},A_{2}\right \} +\left \{ A_{3},A_{4}\right \}
=\left
\{ A_{1}+A_{3},A_{2}+A_{4}\right \} .$

\item $\left \{ A_{1}+A_{3},A_{2}\right \} =\left \{
A_{3}+A_{1},A_{2}\right
\} ,\left \{ A_{1},A_{2}+A_{4}\right \} =\left \{
A_{1},A_{4}+A_{2}\right \} .$

\item $\left \{ A_{1},A_{2}\right \} \left \{ B_{1},B_{2}\right \}
=\{A_{1}B_{1}+A_{2}^{\#}B_{2},A_{1}^{\#}B_{2}+A_{2}B_{1}\}.$

\item $\left( \left \{ A_{1},A_{2}\right \} \left \{ A_{3},A_{4}\right \}
\right) \left \{ A_{5},A_{6}\right \} =\left \{ A_{1},A_{2}\right \} \left(
\left \{ A_{3},A_{4}\right \} \left \{ A_{5},A_{6}\right \} \right) .$
\end{enumerate}

It follows that $y=\left \{ A_{1},A_{2}\right \} x$ is a linear mapping over
the field of real numbers. $\left \{ A_{1},A_{2}\right \} \in \{ \mathbf{C}%
^{n\times m},\mathbf{C}^{n\times m}\}$ is said to be a zero bimatrix
(denoted by $\mathcal{O}_{n\times m}$) if $\left \{ A_{1},A_{2}\right \} x=0$
for any $x\in \mathbf{C}^{m}.$ $\left \{ A_{1},A_{2}\right \} $ is said to
be a square bimatrix if $n=m$. A square bimatrix $\left \{
A_{1},A_{2}\right
\} $ is an identity bimatrix (denoted by $\mathcal{I}_{n}$%
) if $\left \{ A_{1},A_{2}\right \} x=x$, $\forall x\in \mathbf{C}^{n}.$ It
follows that%
\begin{align}
\left \{ A_{1},A_{2}\right \} & =\mathcal{I}_{n}\Longleftrightarrow \left(
A_{1},A_{2}\right) =\left( I_{n},0_{n\times n}\right) ,  \notag \\
\left \{ A_{1},A_{2}\right \} & =\mathcal{O}_{n\times m}\Longleftrightarrow
\left( A_{1},A_{2}\right) =\left( 0_{n\times m},0_{n\times m}\right) .
\label{eq75}
\end{align}%
By Property 5), the power of the square bimatrix $\left \{
A_{1},A_{2}\right
\} $, denoted by $\left \{ A_{1},A_{2}\right \} ^{i}$
with $i\in \mathbf{Z}^{+},$ can be defined recursively as $\left \{
A_{1},A_{2}\right \} ^{i}=\left \{ A_{1},A_{2}\right \} \left \{
A_{1},A_{2}\right \} ^{i-1}$ with $\left \{ A_{1},A_{2}\right \} ^{0}=%
\mathcal{I}_{n}$. It is easy to show, for any $i,j\in \mathbf{Z}^{+},$ that%
\begin{equation*}
\left \{ A_{1},A_{2}\right \} ^{i+j}=\left \{ A_{1},A_{2}\right \} ^{i}\left
\{ A_{1},A_{2}\right \} ^{j}=\left \{ A_{1},A_{2}\right \} ^{j}\left \{
A_{1},A_{2}\right \} ^{i}.
\end{equation*}%
If there exists another bimatrix $\left \{ A_{3},A_{4}\right \} $ such that $%
\left \{ A_{1},A_{2}\right \} \left \{ A_{3},A_{4}\right \} =\left \{
A_{3},A_{4}\right \} \left \{ A_{1},A_{2}\right \} =\mathcal{I}_{n},$ then $%
\left \{ A_{3},A_{4}\right \} $ is said to be the inverse bimatrix of $%
\left
\{ A_{1},A_{2}\right \} $, and is denoted by $\left \{
A_{3},A_{4}\right \} =\left \{ A_{1},A_{2}\right \} ^{-1}.$ For any bimatrix
$\left \{ A_{1},A_{2}\right \} \in \{ \mathbf{C}^{n\times m},\mathbf{C}%
^{n\times m}\},$ denote%
\begin{equation}
\left \{ A_{1},A_{2}\right \} _{\diamond }\triangleq \left[
\begin{array}{cc}
A_{1} & A_{2}^{\#} \\
A_{2} & A_{1}^{\#}%
\end{array}%
\right] \in \mathbf{C}^{2n\times 2m}.  \label{eq88}
\end{equation}%
In the quantum control literature, this is referred to as the so-called
doubled-up matrix associated with the matrix pair $\left(
A_{1},A_{2}^{\#}\right) $ (see \cite{grn10pr,ny17book} and the references
therein). By Property 5), we have $\left \{ A_{1},A_{2}\right \} \left \{
A_{3},A_{4}\right \} =\left \{ X_{1},X_{2}\right \} ,$ with%
\begin{equation*}
\left[
\begin{array}{c}
X_{1} \\
X_{2}%
\end{array}%
\right] =\left \{ A_{1},A_{2}\right \} _{\diamond }\left[
\begin{array}{c}
A_{3} \\
A_{4}%
\end{array}%
\right] ,
\end{equation*}%
by which the following conclusion follows.

\begin{lemma}
\label{lm3}The square bimatrix $\left \{ A_{1},A_{2}\right \} \in \{ \mathbf{%
C}^{n\times n},\mathbf{C}^{n\times n}\}$ is nonsingular if and only if $%
\left \{ A_{1},A_{2}\right \} _{\diamond }$ is nonsingular. Moreover, $%
\left
\{ A_{1},A_{2}\right \} ^{-1}=\left \{ A_{3},A_{4}\right \} $ where%
\begin{equation*}
\left[
\begin{array}{c}
A_{3} \\
A_{4}%
\end{array}%
\right] =\left( \left \{ A_{1},A_{2}\right \} _{\diamond }\right) ^{-1}\left[
\begin{array}{c}
I_{n} \\
0_{n\times n}%
\end{array}%
\right] .
\end{equation*}
\end{lemma}

It can be shown that%
\begin{equation}
\left \{ A_{1},A_{2}\right \} ^{-1}=\left \{
\begin{array}{lc}
\left \{ S_{1}^{-1},-A_{1}^{-\#}A_{2}S_{1}^{-1}\right \} , & \left \vert
A_{1}\right \vert \neq 0, \\
\left \{ -A_{2}^{-1}A_{1}^{\#}S_{2}^{-1},S_{2}^{-1}\right \} , & \left \vert
A_{2}\right \vert \neq 0,%
\end{array}%
\right.  \label{invbimatrix}
\end{equation}%
where $S_{1}=A_{1}-A_{2}^{\#}A_{1}^{-\#}A_{2}$ and $S_{2}=A_{2}^{%
\#}-A_{1}A_{2}^{-1}A_{1}^{\#}$ (notice that, by a Schur complement, $%
S_{i},i=1,2,$ are nonsingular if and only if $\left \{ A_{1},A_{2}\right \}
_{\diamond }$ is nonsingular). However, $\left \{ A_{1},A_{2}\right \} $ can
still be nonsingular if neither $A_{1}$ nor $A_{2}$ is nonsingular, for
example, $A_{1}=\mathrm{diag}\{0,1\}$ and $A_{2}=\mathrm{diag}\{1,0\}$. It
follows that $\left \{ A_{1},0_{n\times n}\right \} $ is nonsingular if and
only if $A_{1}$ is nonsingular, and $\left \{ 0_{n\times n},A_{2}\right \} $
is nonsingular if and only if $A_{2}$ is nonsingular. Moreover,
\begin{equation}
\left \{ A_{1},0_{n\times n}\right \} ^{-1}=\left \{ A_{1}^{-1},0_{n\times
n}\right \} ,\; \left \{ 0_{n\times n},A_{2}\right \} ^{-1}=\{0_{n\times
n},A_{2}^{-\#}\}.  \label{inv2}
\end{equation}

To explain the matrix $\left \{ A_{1},A_{2}\right \} _{\diamond }$
associated with the bimatrix $\left \{ A_{1},A_{2}\right \} ,$ we define a
mapping%
\begin{equation}
\mathbf{C}^{m}\ni x\mapsto \breve{x}=\frac{1}{\sqrt{2}}\left[
\begin{array}{c}
x \\
x^{\#}%
\end{array}%
\right] \in \mathbf{C}^{2m},  \label{eq3a}
\end{equation}%
which is a linear mapping over the field of real numbers.  Notice that $%
\breve{x}$ is slightly different from the doubled-up column vector of $x$ in
the quantum control literature \cite{grn10pr}, where $x$ is an operator on Hilbert spaces. It follows that $\left \Vert
\breve{x}\right \Vert =\left \Vert x\right \Vert ,\forall x\in \mathbf{C}%
^{m}.$ Notice that $x\ $and $\breve{x}$ are not one-to-one, since not for
any $y\in \mathbf{C}^{2m}$ there exists an $x$ such that $y=\breve{x}.$
Taking $\breve{\cdot}$ on both sides of (\ref{opeator}) gives
\begin{equation}
\breve{y}=\left \{ A_{1},A_{2}\right \} _{\diamond }\breve{x},  \label{eq3}
\end{equation}%
in which $\left \{ A_{1},A_{2}\right \} _{\diamond }$ is involved. Hereafter
$\left \{ A_{1},A_{2}\right \} _{\diamond }$ is referred to as the
complex-lifting of $\left \{ A_{1},A_{2}\right \} $.

Denote the real-representation mapping by%
\begin{equation}
\mathbf{C}^{m}\ni x\mapsto \vec{x}=\left[
\begin{array}{c}
\mathrm{Re}\left( x\right) \\
\mathrm{Im}\left( x\right)%
\end{array}%
\right] \in \mathbf{R}^{2m}.  \label{eq90}
\end{equation}%
It follows that $\vec{\cdot}$ is again a linear mapping over the field of
real numbers, and $\left \Vert \vec{x}\right \Vert =\left \Vert
x\right
\Vert ,\forall x\in \mathbf{C}^{m}$.

\begin{lemma}
\label{lm10}The linear mapping (\ref{opeator}) can be expressed by%
\begin{equation}
\vec{y}=\left \{ A_{1},A_{2}\right \} _{\circ}\vec{x},  \label{eq2}
\end{equation}
where $\left \{ A_{1},A_{2}\right \} _{\circ}\in \mathbf{R}^{2n\times2m}$ is
defined by%
\begin{equation}
\left \{ A_{1},A_{2}\right \} _{\circ}=\left[
\begin{array}{cc}
\mathrm{Re}\left( A_{1}+A_{2}\right) & -\mathrm{Im}\left( A_{1}+A_{2}\right)
\\
\mathrm{Im}\left( A_{1}-A_{2}\right) & \mathrm{Re}\left( A_{1}-A_{2}\right)%
\end{array}
\right] .  \label{eq11}
\end{equation}
\end{lemma}

The matrix $\left \{ A_{1},A_{2}\right \} _{\circ }$ will be referred to as
the real-representation of $\left \{ A_{1},A_{2}\right \} $. To see the
relationship between $\breve{\cdot}$ and $\vec{\cdot},$ we write
\begin{equation}
\breve{x}=\frac{1}{\sqrt{2}}\left[
\begin{array}{cc}
I_{n} & \mathrm{j}I_{n} \\
I_{n} & -\mathrm{j}I_{n}%
\end{array}%
\right] \left[
\begin{array}{c}
\mathrm{Re}\left( x\right) \\
\mathrm{Im}\left( x\right)%
\end{array}%
\right] \triangleq H_{n}\vec{x},  \label{eq77}
\end{equation}%
where $H_{n}^{-1}=H_{n}^{\mathrm{H}},$ namely, $H_{n}$ is a unitary
matrix.$\ $Direct manipulation shows that $H_{n}$ also satisfies%
\begin{equation}
H_{n}^{\#}H_{n}^{\mathrm{H}}=H_{n}H_{n}^{\mathrm{T}}=E_{n}=\left[
\begin{array}{cc}
0 & I_{n} \\
I_{n} & 0%
\end{array}%
\right] .  \label{eq14b}
\end{equation}%
The relationship (\ref{eq77}) also links $\left \{ A_{1},A_{2}\right \}
_{\circ }$ and $\left \{ A_{1},A_{2}\right \} _{\diamond },$ as shown below.

\begin{lemma}
\label{lm7}Let $H_{n}$ and $E_{n}$ be given by (\ref{eq77}) and (\ref{eq14b}%
). Then, for any $\left \{ A_{1},A_{2}\right \} \in \{ \mathbf{C}^{n\times
m},\mathbf{C}^{n\times m}\},$%
\begin{align}
\left \{ A_{1},A_{2}\right \} _{\circ} & =H_{n}^{\mathrm{H}}\left \{
A_{1},A_{2}\right \} _{\diamond}H_{m}.  \label{eq14} \\
E_{m}\left( \left \{ A_{1},A_{2}\right \} _{\diamond}\right) ^{\mathrm{T}} &
=\left( \left \{ A_{1},A_{2}\right \} _{\diamond}\right) ^{\mathrm{H}}E_{n}.
\label{eq14a}
\end{align}
\end{lemma}

The following lemma follows from (\ref{eq11}) (or (\ref{eq88}) and (\ref%
{eq14})) and is very useful in this paper.

\begin{lemma}
\label{lm0}For any real matrix $A\in \mathbf{R}^{2n\times 2m},$ there exists
a unique bimatrix $\left \{ A_{1},A_{2}\right \} \in \{ \mathbf{C}^{n\times
m},\mathbf{C}^{n\times m}\}$ such that $\left \{ A_{1},A_{2}\right \}
_{\circ }=A. $ Moreover, $A_{1}$ and $A_{2}$ are given by
\begin{equation}
\left[
\begin{array}{c}
A_{1} \\
A_{2}%
\end{array}%
\right] =H_{n}AH_{m}^{\mathrm{H}}\left[
\begin{array}{c}
I_{m} \\
0_{m\times m}%
\end{array}%
\right] .  \label{eq99}
\end{equation}
\end{lemma}

The next lemma provides further properties of $\left \{ A_{1},A_{2}\right \}
_{\circ}$ and $\left \{ A_{1},A_{2}\right \} _{\diamond}.$

\begin{lemma}
\label{lm2}Let $\{A_{1},A_{2}\} \in \{ \mathbf{C}^{n\times m},\mathbf{C}%
^{n\times m}\}$ and $\{A_{3},A_{4}\} \in \{ \mathbf{C}^{p\times q},\mathbf{C}%
^{p\times q}\}$ be two given bimatrices. Then%
\begin{align*}
\left( \left \{ A_{1},A_{2}\right \} +\left \{ A_{3},A_{4}\right \} \right)
_{\circ} & =\left \{ A_{1},A_{2}\right \} _{\circ}+\left \{
A_{3},A_{4}\right \} _{\circ}, \\
\left( \left \{ A_{1},A_{2}\right \} +\left \{ A_{3},A_{4}\right \} \right)
_{\diamond} & =\left \{ A_{1},A_{2}\right \} _{\diamond}+\left \{
A_{3},A_{4}\right \} _{\diamond}, \\
\left( \left \{ A_{1},A_{2}\right \} \left \{ A_{3},A_{4}\right \} \right)
_{\circ} & =\left \{ A_{1},A_{2}\right \} _{\circ}\left \{ A_{3},A_{4}\right
\} _{\circ}, \\
\left( \left \{ A_{1},A_{2}\right \} \left \{ A_{3},A_{4}\right \} \right)
_{\diamond} & =\left \{ A_{1},A_{2}\right \} _{\diamond}\left \{
A_{3},A_{4}\right \} _{\diamond},
\end{align*}
where the dimensions $n,m,p,$ and $q$ are assumed to be consistent. Moreover,%
\begin{align}
\left \{ A_{1},A_{2}\right \} _{\circ} & =\left \{ A_{3},A_{4}\right \}
_{\circ}\Longleftrightarrow \left \{ A_{1},A_{2}\right \} =\left \{
A_{3},A_{4}\right \} ,  \label{eq55} \\
\left \{ A_{1},A_{2}\right \} _{\diamond} & =\left \{ A_{3},A_{4}\right \}
_{\diamond}\Longleftrightarrow \left \{ A_{1},A_{2}\right \} =\left \{
A_{3},A_{4}\right \} .  \notag
\end{align}
\end{lemma}

Combining Lemmas \ref{lm3} and \ref{lm2} gives $(\left \{
A_{1},A_{2}\right
\} ^{i})_{\circ }=(\left \{ A_{1},A_{2}\right \} _{\circ
})^{i},(\left \{
A_{1},A_{2}\right
\} ^{i})_{\diamond  }=(\left \{ A_{1},A_{2}\right \} _{\diamond
})^{i},i\in \mathbf{Z}^{+}$ and, if $\left \{ A_{1},A_{2}\right \} $ is
nonsingular, then%
\begin{equation}
\left( \left \{ A_{1},A_{2}\right \} ^{-1}\right) _{\circ }=\left( \left \{
A_{1},A_{2}\right \} _{\circ }\right) ^{-1},\left( \left \{
A_{1},A_{2}\right \} ^{-1}\right) _{\diamond }=\left( \left \{
A_{1},A_{2}\right \} _{\diamond }\right) ^{-1}.  \label{inv}
\end{equation}%
We define the conjugate-transpose of the bimatrix $\left \{
A_{1},A_{2}\right \} \in \{ \mathbf{C}^{n\times m},\mathbf{C}^{n\times m}\}$
by
\begin{equation}
\left \{ A_{1},A_{2}\right \} ^{\mathrm{H}}\triangleq \left \{ A_{1}^{%
\mathrm{H}},A_{2}^{\mathrm{T}}\right \} =\left \{ A_{1}^{\mathrm{H}%
},A_{2}^{\# \mathrm{H}}\right \} .  \label{eqtrans}
\end{equation}

\begin{lemma}
\label{lm9}Let $\left \{ A_{1},A_{2}\right \} \in \{ \mathbf{C}^{n\times m},%
\mathbf{C}^{n\times m}\}$ be a given bimatrix. Then%
\begin{equation*}
\left( \left \{ A_{1},A_{2}\right \} ^{\mathrm{H}}\right) _{\circ }=\left(
\left \{ A_{1},A_{2}\right \} _{\circ }\right) ^{\mathrm{H}},\; \left( \left
\{ A_{1},A_{2}\right \} ^{\mathrm{H}}\right) _{\diamond }=\left( \left \{
A_{1},A_{2}\right \} _{\diamond }\right) ^{\mathrm{H}}.
\end{equation*}
\end{lemma}

The square bimatrix $\left \{ P_{1},P_{2}\right \} \in \{ \mathbf{C}%
^{n\times n},\mathbf{C}^{n\times n}\}$ is said to be Hermite if
\begin{equation}
\left \{ P_{1},P_{2}\right \} =\left \{ P_{1},P_{2}\right \} ^{\mathrm{H}%
}=\left \{ P_{1}^{\mathrm{H}},P_{2}^{\mathrm{T}}\right \} .  \label{eq81}
\end{equation}%
Notice that, if $\left \{ P_{1},P_{2}\right \} $ is Hermite, then, for any $%
x\in \mathbf{C}^{n},$
\begin{eqnarray*}
x^{\mathrm{H}}\left \{ P_{1},P_{2}\right \} x &=&x^{\mathrm{H}}\left(
P_{1}x+P_{2}^{\#}x^{\#}\right) \\
&=&x^{\mathrm{H}}\left( P_{1}^{\mathrm{H}}x+P_{2}^{\mathrm{H}}x^{\#}\right)
\\
&=&x^{\mathrm{H}}\left \{ P_{1},P_{2}\right \} ^{\mathrm{H}}x.
\end{eqnarray*}%
However, this does not imply that $x^{\mathrm{H}}\left \{
P_{1},P_{2}\right
\} x$ is a real number (for example, $P_{1}=1$ and $P_{2}=%
\mathrm{j}$) and thus we cannot use $x^{\mathrm{H}}\left \{
P_{1},P_{2}\right \} x>0,\forall x\neq 0$ to define positive definiteness of
$\left \{ P_{1},P_{2}\right \} $. Alternatively, $\left \{
P_{1},P_{2}\right
\} $ is said to be (semi) positive definite (denoted by $%
\left \{ P_{1},P_{2}\right \} >(\geq )0$), if it is an Hermite bimatrix and,
for any $x\in \mathbf{C}^{n},$%
\begin{equation}
\mathrm{Re}\left( x^{\mathrm{H}}\left \{ P_{1},P_{2}\right \} x\right)
>\left( \geq \right) 0,\; \forall x\neq 0.  \label{eq82}
\end{equation}

\begin{lemma}
\label{lm1}The following three statements are equivalent: 1). $\left \{
P_{1},P_{2}\right \} >\left( \geq \right) 0$. 2). $\left \{
P_{1},P_{2}\right \} _{\circ}>\left( \geq \right) 0$. 3). $\left \{
P_{1},P_{2}\right \} _{\diamond}>\left( \geq \right) 0$. Moreover, a matrix $%
P\in \mathbf{R}^{2n\times2n}$ is positive definite if and only if there
exists a unique $\left \{ P_{1},P_{2}\right \} >\left( \geq \right) 0$ such
that $\left \{ P_{1},P_{2}\right \} _{\circ}=P.$
\end{lemma}

We immediately obtain the following corollary of Lemma \ref{lm1}.

\begin{corollary}
\label{coro7}Assume that $Q=Q^{\mathrm{H}}\in \mathbf{C}^{n\times n}.$ Then,
for any $x\in \mathbf{C}^{n},$ $x^{\mathrm{H}}Qx=\vec{x}^{\mathrm{T}%
}\left
\{ Q,0\right \} _{\circ }\vec{x}.$ Moreover, $Q>\left( \geq \right)
0 $ if and only if $\left \{ Q,0\right \} _{\circ }>\left( \geq \right) 0$.
\end{corollary}

\subsection{\label{sec2.2}Alternative representations of the system}

With the bimatrix defined in (\ref{opeator}), the complex-valued linear
system (\ref{sys}) can be written as%
\begin{equation}
\left \{
\begin{array}{rl}
x^{+} & =\left \{ A_{1},A_{2}\right \} x+\left \{ B_{1},B_{2}\right \} u, \\
y & =\left \{ C_{1},C_{2}\right \} x+\left \{ D_{1},D_{2}\right \} u.%
\end{array}%
\right.  \label{sys2}
\end{equation}%
By the real-representation of the bimatrix, system (\ref{sys2}) can also be
expressed by
\begin{equation}
\left \{
\begin{array}{rl}
\vec{x}^{+} & =\left \{ A_{1},A_{2}\right \} _{\circ }\vec{x}+\left \{
B_{1},B_{2}\right \} _{\circ }\vec{u}, \\
\vec{y} & =\left \{ C_{1},C_{2}\right \} _{\circ }\vec{x}+\left \{
D_{1},D_{2}\right \} _{\circ }\vec{u},%
\end{array}%
\right.  \label{realsysm}
\end{equation}%
which is a $2n$ dimensional real-valued linear system, where $\vec{x}\left(
0\right) =\vec{x_{0}}$. Notice that, as $x$ and $\vec{x}$ are one-to-one,$\ $%
systems (\ref{sys}) and (\ref{realsysm}) are equivalent. Hereafter, (\ref%
{realsysm}) is referred to as the real-representation system. A linear
quantum system in (\ref{realsysm}) is referred to as in the quadrature form %
\cite{ny17book}.

\begin{remark}
It follows that a normal $2n$ dimensional real-valued linear system with $2m$
inputs and $2p$ outputs can be equivalently expressed by an $n$ dimensional
complex-valued linear system as (\ref{sys}) with $m$ inputs and $p$ outputs,
namely, dimensions are reduced. This fact has been mentioned in
Introduction, and will be studied in the future.
\end{remark}

In the same way, with the complex-lifting of the bimatrix, we can also
express system (\ref{sys}) by%
\begin{equation}
\left \{
\begin{array}{rl}
\breve{x}^{+} & =\left \{ A_{1},A_{2}\right \} _{\diamond }\breve{x}+\left
\{ B_{1},B_{2}\right \} _{\diamond }\breve{u}, \\
\breve{y} & =\left \{ C_{1},C_{2}\right \} _{\diamond }\breve{x}+\left \{
D_{1},D_{2}\right \} _{\diamond }\breve{u},%
\end{array}%
\right.  \label{complexsys}
\end{equation}%
which is a $2n$ dimensional normal linear system with complex-valued
coefficients, where $\breve{x}\left( 0\right) =\breve{x}_{0}$. However, as $%
x $ and $\breve{x}$ are not one to one, system (\ref{complexsys}) with
\textit{general} coefficients (namely, $\left \{ A_{1},A_{2}\right \}
_{\diamond }$ is replaced by any $2n\times 2n$ complex matrix) may not be
equivalent to system (\ref{sys}). This is easy to understand since $\left \{
A_{1},A_{2}\right \} _{\diamond }$ possesses a special structure as
indicated by (\ref{eq88}). Hereafter, (\ref{complexsys}) is referred to as
the complex-lifting system. A linear quantum system in (\ref{complexsys}) is
referred to as in the doubled-up form \cite{grn10pr,ny17book} (with a
slight difference as $\breve{x}$ here has a different meaning). It is
noticed that coefficient-lifting is a common method in control theory, for
example, for studying periodic coefficient systems \cite{xll18ijss}.

Taking the Laplace (Z) transformation on both sides of (\ref{sys2}) gives%
\begin{equation*}
\left \{
\begin{array}{l}
\left \{ sI_{n}-A_{1},-A_{2}\right \} X=\left \{ B_{1},B_{2}\right \} U, \\
Y=\left \{ C_{1},C_{2}\right \} X+\left \{ D_{1},D_{2}\right \} U,%
\end{array}%
\right.
\end{equation*}%
where and hereafter, unless specifically noted, we should treat the symbol $%
s $ as a \textit{real} parameter (namely, $s^{\#}=s$), and $X,Y,U$
correspond to $x,y,u$ in the frequency domain. Solve the first equation to
get $X=\left \{ sI_{n}-A_{1},-A_{2}\right \} ^{-1}\left \{
B_{1},B_{2}\right
\} U,$ and from the second equation we have $Y=\left \{
G_{1}\left( s\right) ,G_{2}\left( s\right) \right \} U$ where%
\begin{equation*}
\left \{ G_{1}\left( s\right) ,G_{2}\left( s\right) \right \} =\left \{
C_{1},C_{2}\right \} \left \{ sI_{n}-A_{1},-A_{2}\right \} ^{-1}\left \{
B_{1},B_{2}\right \} +\left \{ D_{1},D_{2}\right \} ,
\end{equation*}%
which is referred to as the transfer function bimatrix of system (\ref{sys}%
). Denote the transfer functions of (\ref{realsysm}) and (\ref{complexsys})
by $G_{R}\left( s\right) $ and $G_{C}\left( s\right) ,$ respectively, then
it is easy to see that%
\begin{equation*}
G_{R}\left( s\right) =\left \{ G_{1}\left( s\right) ,G_{2}\left( s\right)
\right \} _{\circ },\;G_{C}\left( s\right) =\left \{ G_{1}\left( s\right)
,G_{2}\left( s\right) \right \} _{\diamond }.
\end{equation*}%
Notice that the first relation implies that $G_{R}\left( s\right) $ is real,
which implies that $s$ should be treated as a \textit{real} parameter. As a
result, such a $G_{C}\left( s\right) $ coincides with the two-side Laplace
transformation in the linear quantum systems theory (see Eq. (10) in \cite{zj13tac}).

We next introduce the eigenvalue of the complex-valued linear system (\ref%
{sys}). For the normal linear system, a scalar $s$ is said to be an
eigenvalue of the system if it is an eigenvalue of $A_{1}.$ However, for the
complex-valued linear system (\ref{sys}), since there are two matrices $%
A_{1} $ and $A_{2}$ involved, we can not define its eigenvalue in the same
way. To this end, we say that a scalar $s$ is said to be an eigenvalue of
system (\ref{sys}) if $s$ is an eigenvalue of its real-representation system
(\ref{realsysm}), or equivalently, an eigenvalue of its complex-lifting
system (\ref{complexsys}). Denote the set of eigenvalues of system (\ref{sys}%
) by $\lambda \left \{ A_{1},A_{2}\right \} ,$ which is symmetric with
respect to the real axis. Thus%
\begin{equation*}
\lambda \left( \left \{ A_{1},A_{2}\right \} \right) =\lambda \left( \left
\{ A_{1},A_{2}\right \} _{\circ }\right) =\lambda \left( \left \{
A_{1},A_{2}\right \} _{\diamond }\right) .
\end{equation*}%
With this definition, we know that the eigenvalue set for the normal linear
system (\ref{normal}) is%
\begin{equation}
\lambda \left( \left \{ A_{1},0\right \} \right) =\lambda \left( \left \{
A_{1},0\right \} _{\diamond }\right) =\lambda \left( A_{1}\right) \cup
\lambda (A_{1}^{\#}),  \label{expanded}
\end{equation}%
which is an expansion of the usual eigenvalue set $\lambda \left(
A_{1}\right) $. Such an expansion is reasonable since it coincides with the
well-known result that system (\ref{normal}) is asymptotically stable if and
only if $A_{1}$ is Hurwitz (Schur). We mention that this definition of
eigenvalue set for the normal linear system (\ref{normal}) (with
complex-valued coefficients) is different from \cite{bodson14cdc,db16tie}
where the eigenvalue set is $\lambda \left( A_{1}\right) ,$ which is
generally not symmetric with the real axis. Similarly, the eigenvalue set
for the antilinear system (\ref{antilinear}) is%
\begin{equation*}
\lambda \left( \left \{ 0,A_{2}\right \} \right) =\lambda \left( \left \{
0,A_{2}\right \} _{\diamond }\right) =\left \{ s:\left \vert \left[
\begin{array}{cc}
sI_{n} & -A_{2}^{\#} \\
-A_{2} & sI_{n}%
\end{array}%
\right] \right \vert =0\right \} .
\end{equation*}

\section{\label{sec3}Analysis of complex-valued linear systems}

\subsection{Solutions of the system equation}

To find solutions to the system equation of (\ref{sys}), we can of course
use the real-representation system (\ref{realsysm}) to get $\vec{x}$ and
then $x.$ However, we are interested in dealing with system (\ref{sys})
directly. To this end, for a bimatrix $\{A_{1},A_{2}\} \in \{ \mathbf{C}%
^{n\times n},\mathbf{C}^{n\times n}\},$ we define the exponent of $\left \{
A_{1},A_{2}\right \} $ by
\begin{equation}
\mathrm{e}^{t\left \{ A_{1},A_{2}\right \} }=\exp \left( t\left \{
A_{1},A_{2}\right \} \right) =\sum \limits_{i=0}^{\infty }\frac{t^{i}\left
\{ A_{1},A_{2}\right \} ^{i}}{i!},t\in \mathbf{R},  \label{exp}
\end{equation}%
which is bounded (and thus well defined) since, for any $x\in \mathbf{C}%
^{n}, $%
\begin{align*}
\left \Vert \exp \left( t\left \{ A_{1},A_{2}\right \} \right) x\right \Vert
& \leq \sum \limits_{i=0}^{\infty }\frac{t^{i}}{i!}\left( \left \Vert
A_{1}\right \Vert +\left \Vert A_{2}^{\#}\right \Vert \right) ^{i}\left
\Vert x\right \Vert \\
& =\exp \left( \left( \left \Vert A_{1}\right \Vert +\left \Vert A_{2}\right
\Vert \right) t\right) \left \Vert x\right \Vert ,
\end{align*}%
where we have used, for any $i\in \mathbf{Z}^{+},$%
\begin{align*}
\left \Vert \left \{ A_{1},A_{2}\right \} ^{i}x\right \Vert & =\left \Vert
\left \{ A_{1},A_{2}\right \} \left \{ A_{1},A_{2}\right \} ^{i-1}x\right
\Vert \\
& \leq \left( \left \Vert A_{1}\right \Vert +\left \Vert A_{2}^{\#}\right
\Vert \right) \left \Vert \left \{ A_{1},A_{2}\right \} ^{i-1}x\right \Vert
\\
& \leq \cdots \leq \left( \left \Vert A_{1}\right \Vert +\left \Vert
A_{2}^{\#}\right \Vert \right) ^{i}\left \Vert x\right \Vert .
\end{align*}

The defined bimatrix exponent possesses similar properties as the normal
matrix exponent. For example:

\begin{itemize}
\item $\exp \left( 0\left \{ A_{1},A_{2}\right \} \right) =\mathcal{I}_{n}.$

\item $\exp \left( \left( t+s\right) \left \{ A_{1},A_{2}\right \} \right)
=\exp \left( t\left \{ A_{1},A_{2}\right \} \right) \exp \left( s\left \{
A_{1},A_{2}\right \} \right) ,\forall t,s\in \mathbf{R}.$

\item $\exp \left( t\left \{ A_{1},A_{2}\right \} \right) $ is nonsingular
with $(\exp \left( t\left \{ A_{1},A_{2}\right \} \right) )^{-1}=\exp \left(
-t\left \{ A_{1},A_{2}\right \} \right) .$

\item $\frac{\mathrm{d}}{\mathrm{d}t}\exp \left( t\left \{
A_{1},A_{2}\right
\} \right) =\left \{ A_{1},A_{2}\right \} \exp \left(
t\left \{ A_{1},A_{2}\right \} \right) .$
\end{itemize}

Equivalent descriptions of
the exponent and power of the bimatrix $\left \{ A_{1},A_{2}\right \} $ can
be obtained as follows.

\begin{lemma}
\label{lm4}Let the pair of matrix functions $\left( \Phi _{1}(t),\Phi
_{2}(t)\right) \in \left( \mathbf{C}^{n\times n},\mathbf{C}^{n\times
n}\right) $ be given by%
\begin{equation}
\left[
\begin{array}{c}
\Phi _{1}(t) \\
\Phi _{2}(t)%
\end{array}%
\right] =\left \{
\begin{array}{ll}
\mathrm{e}^{t\left \{ A_{1},A_{2}\right \} _{\diamond }}\left[
\begin{array}{c}
I_{n} \\
0_{n\times n}%
\end{array}%
\right] , & t\in \mathbf{R}, \\
\left \{ A_{1},A_{2}\right \} _{\diamond }^{t}\left[
\begin{array}{c}
I_{n} \\
0_{n\times n}%
\end{array}%
\right] , & t\in \mathbf{Z},%
\end{array}%
\right.  \label{eq17}
\end{equation}%
where $\left \{ A_{1},A_{2}\right \} $ is assumed to be nonsingular if $t<0.$
Then
\begin{equation}
\left \{ \Phi _{1}(t),\Phi _{2}(t)\right \} =\left \{
\begin{array}{ll}
\mathrm{e}^{t\left \{ A_{1},A_{2}\right \} }, & t\in \mathbf{R}, \\
\left \{ A_{1},A_{2}\right \} ^{t}, & t\in \mathbf{Z}.%
\end{array}%
\right.  \label{eq16}
\end{equation}
\end{lemma}

When the antilinear system (\ref{antilinear}) is concerned, exponent and
power of the bimatrix $\left \{ A_{1},A_{2}\right \} $ can be simplified, as
shown in the next corollary.

\begin{corollary}
\label{coro1}For the antilinear system (\ref{antilinear}), exponent and
power of the bimatrix $\left \{ A_{1},A_{2}\right \} $ can be computed by (%
\ref{eq16}) where
\begin{equation}
\Phi_{1}(t)=\left \{
\begin{array}{ll}
\sum \limits_{i=0}^{\infty}\frac{t^{2i}}{\left( 2i\right) !}\left(
A_{2}^{\#}A_{2}\right) ^{i}, & t\in \mathbf{R}^{+}, \\
\left( A_{2}^{\#}A_{2}\right) ^{\frac{t}{2}}, & t\in \mathbf{Z}^{+},\text{ }t%
\text{ is even,} \\
0, & t\in \mathbf{Z}^{+},\text{ }t\text{ is odd,}%
\end{array}
\right.  \label{eqfif2}
\end{equation}
and%
\begin{equation}
\; \Phi_{2}(t)=\left \{
\begin{array}{ll}
\sum \limits_{i=0}^{\infty}\frac{t^{2i+1}}{\left( 2i+1\right) !}A_{2}\left(
A_{2}^{\#}A_{2}\right) ^{i}, & t\in \mathbf{R}^{+}, \\
0, & t\in \mathbf{Z}^{+},\text{ }t\text{ is even,} \\
A_{2}\left( A_{2}^{\#}A_{2}\right) ^{\frac{t-1}{2}}, & t\in \mathbf{Z}^{+},%
\text{ }t\text{ is odd.}%
\end{array}
\right.  \label{eqf2}
\end{equation}
\end{corollary}

With the exponent and power of the bimatrix $\left \{ A_{1},A_{2}\right \} $%
, solutions to the system equation of (\ref{sys}) can be obtained in
closed-form.

\begin{theorem}
\label{th0}For any $x_{0}\in \mathbf{C}$ and $u(t)\in \mathbf{C}^{m},$ the
system equation of (\ref{sys}) with $x\left( 0\right) =x_{0}$ has the unique
solution given by, for all$\ t\geq0,$%
\begin{align}
x(t)= & \left \{
\begin{array}{l}
\mathrm{e}^{t\left \{ A_{1},A_{2}\right \} }x_{0}+\int_{0}^{t}\mathrm{e}%
^{\left( t-s\right) \left \{ A_{1},A_{2}\right \} }\left \{
B_{1},B_{2}\right \} u\left( s\right) \mathrm{d}s, \\
\left \{ A_{1},A_{2}\right \} ^{t}x_{0}+\sum \limits_{i=0}^{t-1}\left \{
A_{1},A_{2}\right \} ^{t-1-i}\left \{ B_{1},B_{2}\right \} u\left( i\right) ,%
\end{array}
\right.  \label{eq70} \\
= & \left \{ \Phi_{1}(t),\Phi_{2}(t)\right \} x_{0}+\left \{
\begin{array}{l}
\int_{0}^{t}\left \{ \Phi_{1}\left( t-s\right) ,\Phi_{2}\left( t-s\right)
\right \} \left \{ B_{1},B_{2}\right \} u\left( s\right) \mathrm{d}s, \\
\sum \limits_{i=0}^{t-1}\left \{ \Phi_{1}\left( t-1-i\right) ,\Phi_{2}\left(
t-1-i\right) \right \} \left \{ B_{1},B_{2}\right \} u\left( i\right) ,%
\end{array}
\right.  \label{eq78}
\end{align}
where the first row for $t\in \mathbf{R}^{+}$ and the second row for $1\leq
t\in \mathbf{Z}^{+}.$
\end{theorem}

By this theorem it is not
hard to show that this solution coincides with the classical solution to the
normal linear system (\ref{normal}). Applying Theorem \ref{th0} and
Corollary \ref{coro1} on the antilinear system (\ref{antilinear}) gives the
following closed-form solutions:%
\begin{align}
x(t)& =\Phi _{1}(t)x_{0}+\Phi _{2}^{\#}(t)x_{0}^{\#}  \notag \\
& +\left \{
\begin{array}{ll}
\int_{0}^{t}\left( \Phi _{2}^{\#}\left( t-s\right) B_{2}u\left( s\right)
+\Phi _{1}\left( t-s\right) B_{2}^{\#}u^{\#}\left( s\right) \right) \mathrm{d%
}s, & t\in \mathbf{R}^{+}, \\
\sum \limits_{i=0}^{t-1}\left( \Phi _{2}^{\#}\left( t-1-i\right)
B_{2}u\left( i\right) +\Phi _{1}\left( t-1-i\right) B_{2}^{\#}u^{\#}\left(
i\right) \right) , & 1\leq t\in \mathbf{Z}^{+},%
\end{array}%
\right.  \label{eqsolu}
\end{align}%
where $(\Phi _{1}(t),\Phi _{2}(t))$ is given by (\ref{eqfif2})-(\ref{eqf2}).
Solutions (\ref{eqsolu}) seem easier to recognize than the equations
obtained in \cite{wzls15iet} where $t\in \mathbf{R}^{+}$ was studied.

\subsection{Controllability, observability and stability}

If for any initial condition $x_{0}\in \mathbf{C}^{n}$ and any final state $%
x_{1}\in \mathbf{C}^{n}$, there exists a function $u(t)\in $ $\mathbf{C}^{m}$
and a finite number $t_{1}>0$ such that the solution satisfies $x\left(
t_{1}\right) =x_{1},$ then system (\ref{sys}) is said to be controllable (or
$(\left \{ A_{1},A_{2}\right \} ,\left \{ B_{1},B_{2}\right \} )$ is
controllable). In a similar way, we can define the observability of system (%
\ref{sys}). The system (\ref{sys}) is said to be stable if, for any $\delta
>0$ there exists an $\varepsilon $ such that $\left \Vert x_{0}\right \Vert
\leq \varepsilon \Rightarrow \left \Vert x(t)\right \Vert \leq \delta ,\forall
t\geq 0;$ is said to be asymptotically stable if it is stable and, moreover,
for any $\epsilon >0$ there exists a number $\eta $ such that $\left \Vert
x(t)\right \Vert \leq \epsilon ,\forall t\geq \eta .$ Since $x\mapsto \vec{x}$
is one-to one, controllability, observability and stability of system (\ref%
{sys}) are equivalent to that of its real-representation system (\ref%
{realsysm}), which, by (\ref{eq14}), are further equivalent to that of its
complex-lifting system (\ref{complexsys}).

Thus, by the PBH test \cite{kfa69book}, system (\ref{sys}) is controllable
if and only if%
\begin{equation}
\mathrm{rank}\left[
\begin{array}{cccc}
sI_{n}-A_{1} & -A_{2}^{\#} & B_{1} & B_{2}^{\#} \\
-A_{2} & sI_{n}-A_{1}^{\#} & B_{2} & B_{1}^{\#}%
\end{array}%
\right] =2n,\; \forall s\in \mathbf{C},  \label{pbhtest}
\end{equation}%
is observable if and only if%
\begin{equation}
\mathrm{rank}\left[
\begin{array}{cc}
sI_{n}-A_{1} & -A_{2}^{\#} \\
-A_{2} & sI_{n}-A_{1}^{\#} \\
C_{1} & C_{2}^{\#} \\
C_{2} & C_{1}^{\#}%
\end{array}%
\right] =2n,\; \forall s\in \mathbf{C},  \label{phbtest}
\end{equation}%
and is asymptotically stable if and only if $\left \{ A_{1},A_{2}\right \}
_{\diamond }$ is Hurwitz (Schur when $t\in \mathbf{Z}^{+}$).

It is readily to see that the above conclusion coincides with classical
results for the normal linear system (\ref{normal}), namely, $(\left \{
A_{1},0\right \} ,\left \{ B_{1},0\right \} )$ is controllable (observable) if
and only if $\left( A_{1},B_{1}\right) $ is controllable (observable), $%
\left \{ A_{1},0\right \} $ is asymptotically stable if and only if $A_{1}$ is
Hurwitz (Schur when $t\in \mathbf{Z}^{+}$). However, different situation
exists for the antilinear system (\ref{antilinear}), as made clear in the
following corollary.

\begin{corollary}
\label{coro1a}The antilinear system (\ref{antilinear}) is controllable if
and only if%
\begin{equation}
\mathrm{rank}\left[
\begin{array}{ccc}
sI_{n}-A_{2}^{\#}A_{2} & B_{2}^{\#} & A_{2}^{\#}B_{2}%
\end{array}
\right] =n,\; \forall s\in \mathbf{C},  \label{eq18}
\end{equation}
and is observable if and only if%
\begin{equation}
\left[
\begin{array}{c}
sI_{n}-A_{2}^{\#}A_{2} \\
C_{2} \\
C_{2}^{\#}A_{2}%
\end{array}
\right] =n,\; \forall s\in \mathbf{C}.  \label{eq19}
\end{equation}
\end{corollary}

It is very interesting to emphasize that the antilinear system (\ref%
{antilinear}) is controllable (observable) if and only if the associated
normal linear system $(A_{2}^{\#}A_{2},[B_{2}^{\#},A_{2}^{\#}B_{2}])$ is
controllable ($(A_{2}^{\#}A_{2},[C_{2}^{\#};C_{2}^{\#}A_{2}])$ is
observable).

\begin{corollary}
\label{coro2}The antilinear system (\ref{antilinear}) with $t\in \mathbf{Z}%
^{+}$ is asymptotically stable if and only if $A_{2}^{\#}A_{2}$ is Schur,
namely,
\begin{equation}
\rho \left( A_{2}^{\#}A_{2}\right) <1,  \label{spectral}
\end{equation}
and with $t\in \mathbf{R}^{+}$ cannot be asymptotically stable for any $%
A_{2}.$
\end{corollary}

Lyapunov stability theorem for the
complex-valued linear system (\ref{sys}) can also be obtained.

\begin{theorem}
\label{th6}Let the bimatrix $\{C_{1},C_{2}\} \in \{ \mathbf{C}^{p\times n},%
\mathbf{C}^{p\times n}\}$ be such that $(\left \{ A_{1},A_{2}\right \}
,\left \{ C_{1},C_{2}\right \} )$ is observable. Then system (\ref{sys}) is
asymptotically stable if and only if there exists a positive definite
bimatrix $\left \{ P_{1},P_{2}\right \} $ $\in \{ \mathbf{C}^{n\times n},%
\mathbf{C}^{n\times n}\}$ solving the Lyapunov bimatrix equation%
\begin{equation}
-\left \{ C_{1},C_{2}\right \} ^{\mathrm{H}}\left \{ C_{1},C_{2}\right \}
=\left \{
\begin{array}{lc}
\left \{ A_{1},A_{2}\right \} ^{\mathrm{H}}\left \{ P_{1},P_{2}\right \}
+\left \{ P_{1},P_{2}\right \} \left \{ A_{1},A_{2}\right \} , & t\in
\mathbf{R}^{+}, \\
\left \{ A_{1},A_{2}\right \} ^{\mathrm{H}}\left \{ P_{1},P_{2}\right \}
\left \{ A_{1},A_{2}\right \} -\left \{ P_{1},P_{2}\right \} , & t\in
\mathbf{Z}^{+},%
\end{array}
\right.  \label{eq33}
\end{equation}
or, equivalently, the coupled matrix equations%
\begin{equation}
\left \{
\begin{array}{rl}
-\left( C_{1}^{\mathrm{H}}C_{1}+C_{2}^{\mathrm{H}}C_{2}\right) & =\left \{
\begin{array}{l}
A_{1}^{\mathrm{H}}P_{1}+P_{1}A_{1}+A_{2}^{\mathrm{H}}P_{2}+P_{2}^{\mathrm{H}%
}A_{2}, \\
A_{1}^{\mathrm{H}}P_{1}A_{1}+A_{2}^{\mathrm{H}}P_{1}^{\#}A_{2}+A_{2}^{%
\mathrm{H}}P_{2}A_{1}+A_{1}^{\mathrm{H}}P_{2}^{\mathrm{H}}A_{2}-P_{1},%
\end{array}
\right. \\
-\left( C_{1}^{\mathrm{T}}C_{2}+C_{2}^{\mathrm{T}}C_{1}\right) & =\left \{
\begin{array}{l}
A_{1}^{\mathrm{T}}P_{2}+P_{2}A_{1}+A_{2}^{\mathrm{T}}P_{1}+P_{1}^{\#}A_{2},
\\
A_{1}^{\mathrm{T}}P_{1}^{\#}A_{2}+A_{2}^{\mathrm{T}}P_{1}A_{1}+A_{2}^{%
\mathrm{T}}P_{2}^{\mathrm{H}}A_{2}+A_{1}^{\mathrm{T}}P_{2}A_{1}-P_{2},%
\end{array}
\right.%
\end{array}
\right.  \label{eq35}
\end{equation}
where, in each equation, the first row for $t\in \mathbf{R}^{+}$ and the
second row for $t\in \mathbf{Z}^{+}.$
\end{theorem}

Two remarks are given in order regarding Theorem \ref{th6}.

\begin{remark}
The Lyapunov bimatrix equation (\ref{eq33}) can also be written in another
matrix form:
\begin{equation}
-\left \{ C_{1},C_{2}\right \} _{\diamond }^{\mathrm{H}}\left \{
C_{1},C_{2}\right \} _{\diamond }=\left \{
\begin{array}{cc}
\left \{ A_{1},A_{2}\right \} _{\diamond }^{\mathrm{H}}P_{C}+P_{C}\left \{
A_{1},A_{2}\right \} _{\diamond }, & t\in \mathbf{R}^{+}, \\
\left \{ A_{1},A_{2}\right \} _{\diamond }^{\mathrm{H}}P_{C}\left \{
A_{1},A_{2}\right \} _{\diamond }-P_{C}, & t\in \mathbf{Z}^{+},%
\end{array}%
\right.  \label{eq52}
\end{equation}%
which corresponds to the complex-lifting system (\ref{complexsys}). This
equation can be obtained from (\ref{eq33}) by using (\ref{eq14}). Moreover,
solutions to (\ref{eq52}) and (\ref{eq33}) are related with $P_{C}=\left \{
P_{1},P_{2}\right \} _{\diamond }.$
\end{remark}

\begin{remark}
Assume that the Lyapunov bimatrix equation (\ref{eq33}) has a solution $%
\left \{ P_{1},P_{2}\right \} >0.$ Consider a Lyapunov function
\begin{equation}
V\left( x\right) =\mathrm{Re}\left( x^{\mathrm{H}}\left \{ P_{1},P_{2}\right
\} x\right) =\vec{x}^{\mathrm{T}}\left \{ P_{1},P_{2}\right \} _{\circ }\vec{%
x},  \label{lya}
\end{equation}%
which is positive definite. Then, by using (%
\ref{eq33}), the time-derivative (time-shift) of $V\left(
x\right) $ along system (\ref{sys}) (or system (\ref{realsysm})) is given by%
\begin{align*}
V^{+}\left( x\right) & =\left \{
\begin{array}{cc}
\mathrm{Re}\left( x^{\mathrm{H}}\left( \left \{ A_{1},A_{2}\right \} ^{%
\mathrm{H}}\left \{ P_{1},P_{2}\right \} +\left \{ P_{1},P_{2}\right \}
\left \{ A_{1},A_{2}\right \} \right) x\right) , & t\in \mathbf{R}^{+}, \\
\mathrm{Re}\left( x^{\mathrm{H}}\left( \left \{ A_{1},A_{2}\right \} ^{%
\mathrm{H}}\left \{ P_{1},P_{2}\right \} \left \{ A_{1},A_{2}\right \}
-\left \{ P_{1},P_{2}\right \} \right) x\right) , & t\in \mathbf{Z}^{+},%
\end{array}%
\right. \\
& =-\mathrm{Re}\left( x^{\mathrm{H}}\left \{ C_{1},C_{2}\right \} ^{\mathrm{H%
}}\left \{ C_{1},C_{2}\right \} x\right) \\
& =-\vec{x}^{\mathrm{T}}\left( \left \{ C_{1},C_{2}\right \} _{\circ }^{%
\mathrm{H}}\left \{ C_{1},C_{2}\right \} _{\circ }\right) \vec{x}.
\end{align*}%
Since $(\left \{ A_{1},A_{2}\right \} _{\circ },\left \{
C_{1},C_{2}\right
\} _{\circ })$ is observable, it follows from the
well-known Lyapunov stability theorem \cite{rugh96book} that $\left \{
A_{1},A_{2}\right \} _{\circ }$ is asymptotically stable. The converse can
be shown similarly.
\end{remark}

It is not hard to show that Theorem \ref{th6} reduces to the well-known
Lyapunov stability theorem when it is applied on the normal linear system (%
\ref{normal}), namely, $A_{2}=0$, $C_{2}=0$ and $\left( A_{1},C_{1}\right) $
is observable. In fact, the coupled matrix equations in (\ref{eq35}) become
\begin{equation*}
\left \{
\begin{array}{rl}
-C_{1}^{\mathrm{H}}C_{1} & =\left \{
\begin{array}{ll}
A_{1}^{\mathrm{H}}P_{1}+P_{1}A_{1} & t\in \mathbf{R}^{+}, \\
A_{1}^{\mathrm{H}}P_{1}A_{1}-P_{1}, & t\in \mathbf{Z}^{+},%
\end{array}
\right. \\
0 & =\left \{
\begin{array}{ll}
A_{1}^{\mathrm{T}}P_{2}+P_{2}A_{1}, & t\in \mathbf{R}^{+}, \\
A_{1}^{\mathrm{T}}P_{2}A_{1}-P_{2}, & t\in \mathbf{Z},%
\end{array}
\right.%
\end{array}
\right.
\end{equation*}
which have a unique solution $\left \{ P_{1},P_{2}\right \} =\left \{ P_{1}^{%
\mathrm{H}},0\right \} >0$ if and only if $A_{1}$ is asymptotically stable.

\begin{corollary}
\label{coro3}Let $C_{2}\in \mathbf{C}^{p\times n}$ be any matrix such that
system (\ref{antilinear}) is observable and $C_{N}\in \mathbf{C}^{q\times n}$
be such that system $(A_{2}^{\#}A_{2},C_{N})$ is observable. Then the
following statements are equivalent:

\begin{enumerate}
\item The antilinear system (\ref{antilinear}) with $t\in \mathbf{Z}^{+}$ is
asymptotically stable.

\item There exists a unique positive definite matrix $P_{1}\in \mathbf{C}%
^{n\times n}$ such that%
\begin{equation}
A_{2}^{\mathrm{H}}P_{1}^{\#}A_{2}-P_{1}=-C_{2}^{\mathrm{H}}C_{2}.
\label{lya1}
\end{equation}

\item There exists a unique positive definite matrix $P_{N}\in \mathbf{C}%
^{n\times n}$ such that%
\begin{equation}
\left( A_{2}^{\#}A_{2}\right) ^{\mathrm{H}}P_{N}\left(
A_{2}^{\#}A_{2}\right) -P_{N}=-C_{N}^{\mathrm{H}}C_{N}.  \label{lya2}
\end{equation}
\end{enumerate}

Moreover, if system (\ref{antilinear}) is asymptotically stable, $\left(
C_{2},P_{1}\right) $ satisfies (\ref{lya1}) if and only if $\left(
C_{N},P_{N}\right) $ defined below satisfies (\ref{lya2}):
\begin{equation}
\left( C_{N},P_{N}\right) =\left( \left[
\begin{array}{c}
C_{2} \\
C_{2}^{\#}A_{2}%
\end{array}
\right] ,P_{1}\right) .  \label{eq73}
\end{equation}
\end{corollary}

Items 1) and 2) generalize a result in \cite{wdl13aucc} where $C_{2}^{%
\mathrm{H}}C_{2}$ is assumed to be positive definite. We point out that our
proof is also different from \cite{wdl13aucc} where a Lyapunov function
approach was constructed, while our proof is algebraic.

Of course, there is no need to investigate Lyapunov stability theorem for
system (\ref{antilinear}) with $t\in \mathbf{R}^{+}$ in view of Corollary %
\ref{coro2}.

\section{\label{sec4}Design of complex-valued linear systems}

\subsection{\label{sec4.1}Eigenvalue assignment and stabilization}

Similarly to normal linear systems, we may consider eigenvalue assignment
for the complex-valued linear system (\ref{sys}) by the full state feedback (%
\ref{eqfeedback}) where $\left \{ K_{1},K_{2}\right \} $\ $\in \{ \mathbf{C}%
^{m\times n},\mathbf{C}^{m\times n}\}$ is the feedback gain bimatrix. The
resulting closed-loop system is%
\begin{equation}
x^{+}=\left( \left \{ A_{1},A_{2}\right \} +\left \{ B_{1},B_{2}\right \}
\left \{ K_{1},K_{2}\right \} \right) x.  \label{closed1}
\end{equation}%
Then (\ref{eqfeedback}) is said to be an eigenvalue assignment controller if
(\ref{closed1}) possesses the desired eigenvalue set $\Gamma $ that is
symmetric with respect to the real axis. In this case $\left \{
K_{1},K_{2}\right \} $ is said to be the eigenvalue assignment gain bimatrix.

\begin{theorem}
\label{th5}For any given $\Gamma $ that is symmetric with respect to the
real axis, there exists an eigenvalue assignment controller for system (\ref%
{sys}) if and only if it is controllable. Moreover, $\left \{
K_{1},K_{2}\right \} $ is an eigenvalue assignment gain bimatrix for system (%
\ref{sys}) if and only if $K\in \mathbf{R}^{2m\times 2n}$ is an eigenvalue
assignment gain for the real-valued normal linear system $(\left \{
A_{1},A_{2}\right \} _{\circ },\left \{ B_{1},B_{2}\right \} _{\circ }),$
where%
\begin{equation}
\left[
\begin{array}{c}
K_{1} \\
K_{2}%
\end{array}%
\right] =H_{m}KH_{n}^{\mathrm{H}}\left[
\begin{array}{c}
I_{n} \\
0_{n\times n}%
\end{array}%
\right] .  \label{eqk1k2}
\end{equation}
\end{theorem}

By this theorem, to solve the eigenvalue assignment problem for system (\ref%
{sys}), we first use any standard eigenvalue assignment method (see, for
example, \cite{duan93tac, lt97ijss, wf06tspl, zl16ijss}) to compute the
real-valued feedback gain $K$. Then the eigenvalue
assignment gain bimatrix $\left \{ K_{1},K_{2}\right \} $ can be computed
according to (\ref{eqk1k2}). Finally, the full state feedback controller can
be constructed according to (\ref{eqfeedback}). A systematic design of
eigenvalue assignment for complex-valued linear system by using generalized
Sylvester bimatrix equations (see the normal cases studied in %
\cite{duan93tac, zd06scl}, and \cite{zdl09iet}) will be reported elsewhere.

\begin{remark}
\label{rm1}It follows from Theorem \ref{th5} that we can use the
real-representation system (\ref{realsysm}) to design eigenvalue assignment
gain bimatrix $\left \{ K_{1},K_{2}\right \} .$ However, we can not use the
complex-lifting system (\ref{complexsys}) instead, though controllabilities
of these systems are equivalent. This is because, if $K$ is designed such
that $\left \{ A_{1},A_{2}\right \} _{\diamond }+\left \{
B_{1},B_{2}\right
\} _{\diamond }K$ possesses the prescribed eigenvalue set
$\Gamma $, $K$ generally has no evident structure such that there exists a
bimatrix $\left
\{ K_{1},K_{2}\right \} $ satisfying
\begin{equation}
\left \{ K_{1},K_{2}\right \} _{\diamond }=\left[
\begin{array}{cc}
K_{1} & K_{2}^{\#} \\
K_{2} & K_{1}^{\#}%
\end{array}%
\right] =K.  \label{consistent}
\end{equation}%
\end{remark}

We next discuss the stabilization problem for system (\ref{sys}). The full
state feedback controller (\ref{eqfeedback}) is said to be a stabilizing
controller for this system if the closed-loop system (\ref{closed1}) is
asymptotically stable. In this case system (\ref{sys}) is said to be
stabilizable, and $\left \{ K_{1},K_{2}\right \} $ is said to be a
stabilizing gain bimatrix. Similar to the case of controllability, we can
see that system (\ref{sys}) is stabilizable if and only if system (\ref%
{realsysm}) is, which is further equivalent to the stabilizability of system
(\ref{complexsys}), namely,%
\begin{equation*}
\mathrm{rank}\left[
\begin{array}{cccc}
sI_{n}-A_{1} & -A_{2}^{\#} & B_{1} & B_{2}^{\#} \\
-A_{2} & sI_{n}-A_{1}^{\#} & B_{2} & B_{1}^{\#}%
\end{array}%
\right] =2n,\; \forall s\in \mathbf{D},
\end{equation*}%
where $\mathbf{D}=\{s:\mathrm{Re}\{s\} \geq 0\}$ if $t\in \mathbf{R}^{+}$
and $\mathbf{D}=\{s:\left \vert s\right \vert \geq 1\}$ if $t\in \mathbf{Z}%
^{+}.$ Moreover, $K\in \mathbf{R}^{2m\times 2n}$ is a stabilizing gain for
the normal linear system (\ref{realsysm}) if and only if $\left \{
K_{1},K_{2}\right \} $ defined in (\ref{eqk1k2}) is a stabilizing gain
bimatrix for system (\ref{sys}).

For the real-valued normal linear system (\ref{normal}) with the normal
linear feedback (\ref{normalfeedback}), it is known that the degree of
freedom in the design is proportional to the number of inputs $m$ %
\cite{duan93tac, zd06scl}. However, for the complex-valued linear system (%
\ref{spectral}) with the full-state feedback (\ref{eqfeedback}), as the
real-representation system (\ref{realsysm}) has $2m$ inputs, the degree of
freedom in the design has been doubled. This fact indicates that the full
state feedback (\ref{eqfeedback}) can introduce more freedom in the design.
An example has been given in Introduction to illustrate this point.

We finally check the stabilizability of the antilinear system (\ref%
{antilinear}).

\begin{corollary}
\label{coro6}The antilinear system (\ref{antilinear}) with $t\in \mathbf{R}%
^{+}$ is stabilizable if and only if it is controllable, namely, (\ref{eq18}%
) is satisfied, and with $t\in \mathbf{Z}^{+}$ is stabilizable if and only if%
\begin{equation}
\mathrm{rank}\left[
\begin{array}{ccc}
\lambda I_{n}-A_{2}A_{2}^{\#} & B_{2} & A_{2}B_{2}^{\#}%
\end{array}%
\right] =n,\; \forall \lambda \in \{s:\left \vert s\right \vert \geq 1\}.
\label{eq18a}
\end{equation}
\end{corollary}

We mention that,
similarly to system (\ref{sys}), full state feedback (\ref{eqfeedback}) is
generally necessary for stabilizing (or achieving eigenvalue assignment) the
antilinear system (\ref{antilinear}). This is particularly clear when $t\in
\mathbf{R}^{+}$ since the closed-loop system%
\begin{equation}
x^{+}=\left( A_{2}^{\#}+B_{2}^{\#}K_{1}^{\#}\right) x^{\#},  \label{eq43}
\end{equation}%
by the normal state feedback (\ref{normalfeedback}) can not be
asymptotically stable according to Corollary \ref{coro2}.

\subsection{Linear quadratic regulation}

Consider the real-valued index functional%
\begin{equation}
J\left( u\right) =\left \{
\begin{array}{lc}
\int_{0}^{\infty }\mathrm{Re}\left( x^{\mathrm{H}}(t)\left \{
Q_{1},Q_{2}\right \} x(t)+u^{\mathrm{H}}(t)\left \{ R_{1},R_{2}\right \}
u(t)\right) \mathrm{d}t, & t\in \mathbf{R}^{+}, \\
\sum \limits_{t=0}^{\infty }\mathrm{Re}\left( x^{\mathrm{H}}(t)\left \{
Q_{1},Q_{2}\right \} x(t)+u^{\mathrm{H}}(t)\left \{ R_{1},R_{2}\right \}
u(t)\right) , & t\in \mathbf{Z}^{+},%
\end{array}%
\right.  \label{eqj}
\end{equation}%
where $\left \{ Q_{1},Q_{2}\right \} \in \left \{ \mathbf{C}^{n\times n},%
\mathbf{C}^{n\times n}\right \} $ and $\left \{ R_{1},R_{2}\right \} \in
\left \{ \mathbf{C}^{m\times m},\mathbf{C}^{m\times m}\right \} $ are given
positive definite weighting bimatrices ($\left \{ Q_{1},Q_{2}\right \} $ can
be semi-positive definite, however, we assume $\left \{ Q_{1},Q_{2}\right \}
>0 $ for simplicity). The linear quadratic regulation (LQR) problem is
stated as finding an optimal controller $u^{\ast }$ for system (\ref{sys})
such that $J\left( u\right) $ is minimized, denoted by $J_{\min }\left(
u^{\ast }\right) $.

\begin{theorem}
\label{th3}Assume that the complex-valued linear system (\ref{sys}) is
stabilizable. Then there is a unique bimatrix $\left \{ P_{1},P_{2}\right \}
>0 $ to the following bimatrix ARE%
\begin{eqnarray}
&-\left \{ Q_{1},Q_{2}\right \} =&\left \{ A_{1},A_{2}\right \} ^{\mathrm{H}%
}\left \{ P_{1},P_{2}\right \} +\left \{ P_{1},P_{2}\right \} \left \{
A_{1},A_{2}\right \}  \notag \\
&&-\left \{ P_{1},P_{2}\right \} \left \{ B_{1},B_{2}\right \} \left \{
R_{1},R_{2}\right \} ^{-1}\left \{ B_{1},B_{2}\right \} ^{\mathrm{H}}\left
\{ P_{1},P_{2}\right \} ,  \label{eqareA}
\end{eqnarray}%
when $t\in \mathbf{R}^{+}$ and the bimatrix ARE%
\begin{align}
-\left \{ Q_{1},Q_{2}\right \}& =\left \{ A_{1},A_{2}\right \} ^{\mathrm{H}%
}\left \{ P_{1},P_{2}\right \} \left \{ A_{1},A_{2}\right \} -\left \{
P_{1},P_{2}\right \}  \notag \\
&-\left \{ A_{1},A_{2}\right \} ^{\mathrm{H}}\left \{ P_{1},P_{2}\right \}
\left \{ B_{1},B_{2}\right \} \left \{ S_{1},S_{2}\right \} ^{-1}\left \{
B_{1},B_{2}\right \} ^{\mathrm{H}}\left \{ P_{1},P_{2}\right \} \left \{
A_{1},A_{2}\right \} ,  \label{eqareB}
\end{align}%
when $t\in \mathbf{Z}^{+},$ where $\left \{ S_{1},S_{2}\right \} =\left \{
R_{1},R_{2}\right \} +\left \{ B_{1},B_{2}\right \} ^{\mathrm{H}}\left \{
P_{1},P_{2}\right \} \left \{ B_{1},B_{2}\right \} .$ Moreover, the optimal
control is the full state feedback
\begin{equation}
u^{\ast }=\left \{ K_{1}^{\ast },K_{2}^{\ast }\right \} x,
\label{optcontrol}
\end{equation}%
where $\left \{ K_{1}^{\ast },K_{2}^{\ast }\right \} $ is the optimal
feedback gain bimatrix determined by%
\begin{equation}
\left \{ K_{1}^{\ast },K_{2}^{\ast }\right \} =-\left \{
\begin{array}{ll}
\left \{ R_{1},R_{2}\right \} ^{-1}\left \{ B_{1},B_{2}\right \} ^{\mathrm{H}%
}\left \{ P_{1},P_{2}\right \} , & t\in \mathbf{R}^{+}, \\
\left \{ S_{1},S_{2}\right \} ^{-1}\  \left \{ B_{1},B_{2}\right \} ^{\mathrm{%
H}}\left \{ P_{1},P_{2}\right \} \left \{ A_{1},A_{2}\right \} , & t\in
\mathbf{Z}^{+},%
\end{array}%
\right.  \label{eqgain}
\end{equation}%
and the closed-loop system is asymptotically stable with
\begin{equation}
J_{\min }\left( u^{\ast }\right) =\mathrm{Re}\left( x_{0}^{\mathrm{H}}\left
\{ P_{1},P_{2}\right \} x_{0}\right) .  \label{minj}
\end{equation}
\end{theorem}

Several remarks are given in order regarding Theorem \ref{th3}.

\begin{remark}
Similarly to Theorem \ref{th6}, the bimatrix ARE (\ref{eqareA}) and (\ref%
{eqareB}) may also be written as coupled matrix equations. When $t\in
\mathbf{R}^{+}$ and $R_{2}=0,$ this is easy since (\ref{eqareA}) is
equivalent%
\begin{equation}
\left \{
\begin{array}{rl}
-Q_{1}= & A_{1}^{\mathrm{H}}P_{1}+P_{1}A_{1}+A_{2}^{\mathrm{H}%
}P_{2}+P_{2}^{\#}A_{2} \\
& -\left( P_{1}B_{1}+P_{2}^{\#}B_{2}\right) R_{1}^{-1}\left( B_{1}^{\mathrm{H%
}}P_{1}+B_{2}^{\mathrm{H}}P_{2}\right) \\
& -\left( P_{1}B_{2}^{\#}+P_{2}^{\#}B_{1}^{\#}\right) R_{1}^{-\#}\left(
B_{1}^{\mathrm{T}}P_{2}+B_{2}^{\mathrm{T}}P_{1}\right) , \\
-Q_{2}= & A_{1}^{\mathrm{T}}P_{2}+P_{2}A_{1}+A_{2}^{\mathrm{T}%
}P_{1}+P_{1}^{\#}A_{2} \\
& -\left( P_{1}^{\#}B_{2}+P_{2}B_{1}\right) R_{1}^{-1}\left( B_{1}^{\mathrm{H%
}}P_{1}+B_{2}^{\mathrm{H}}P_{2}\right) \\
& -\left( P_{1}^{\#}B_{1}^{\#}+P_{2}B_{2}^{\#}\right) R_{1}^{-\#}\left(
B_{1}^{\mathrm{T}}P_{2}+B_{2}^{\mathrm{T}}P_{1}\right) .%
\end{array}%
\right.  \label{eqare1}
\end{equation}%
However, when $t\in \mathbf{Z}^{+},$ the results are less interesting as the
inverse of $\left \{ S_{1},S_{2}\right \} =\left \{ R_{1},R_{2}\right \}
+\left \{ B_{1},B_{2}\right \} ^{\mathrm{H}}\left \{ P_{1},P_{2}\right \}
\left \{ B_{1},B_{2}\right \} $ is much more involved (see (\ref{invbimatrix}%
)).
\end{remark}

\begin{remark}
\label{rm2}It is easy to show, by using equation (\ref{eq14}) in Lemma \ref%
{lm7}, that the ARE (\ref{eqareA}) and (\ref{eqareB}) are also equivalent to%
\begin{equation}
-Q_{C}=\left \{
\begin{array}{ll}
A_{C}^{\mathrm{H}}P_{C}+P_{C}A_{C}-P_{C}B_{C}R_{C}^{-1}B_{C}^{\mathrm{H}%
}P_{C}, & t\in \mathbf{R}^{+}, \\
A_{C}^{\mathrm{H}}P_{C}A_{C}-P_{C}-A_{C}^{\mathrm{H}%
}P_{C}B_{C}S_{C}^{-1}B_{C}^{\mathrm{H}}P_{C}A_{C}, & t\in \mathbf{Z}^{+},%
\end{array}%
\right.  \label{are1}
\end{equation}%
according to $P_{C}=H_{n}P_{R}H_{n}^{\mathrm{H}}=\left \{
P_{1},P_{2}\right
\} _{\diamond },$ where $A_{C}=\left \{
A_{1},A_{2}\right
\} _{\diamond },B_{C}=\left \{ B_{1},B_{2}\right \}
_{\diamond },Q_{C}=\left
\{ Q_{1},Q_{2}\right \} _{\diamond }\in \mathbf{C}%
^{2n\times 2n},$ $R_{C}=\left
\{ R_{1},R_{2}\right \} _{\diamond }\in
\mathbf{C}^{2m\times 2m} $ and $S_{C}=R_{C}+B_{C}^{\mathrm{H}}P_{C}B_{C}.$
Hence these AREs (\ref{eqareA})--(\ref{eqareB}), and (\ref{are1}%
) are equivalent.
\end{remark}

\begin{remark}
It follows that the optimal feedback gain bimatrix in (\ref%
{optcontrol}) can be also computed as%
\begin{equation}
\left[
\begin{array}{c}
K_{1}^{\ast } \\
K_{2}^{\ast }%
\end{array}%
\right] =K_{C}\left[
\begin{array}{c}
I_{n} \\
0_{n\times n}%
\end{array}%
\right] ,  \label{optk1k2}
\end{equation}%
where%
\begin{equation}
K_{C}=\left \{
\begin{array}{ll}
-R_{C}^{-1}B_{C}^{\mathrm{H}}P_{C}, & t\in \mathbf{R}^{+}, \\
-S_{C}^{-1}B_{C}^{\mathrm{H}}P_{C}A_{C}, & t\in \mathbf{Z}^{+}.%
\end{array}%
\right.  \label{eqkc}
\end{equation}%
Moreover, the the minimal value of $J\left( u\right) $ is given by
\begin{equation}
J_{\min }\left( u^{\ast }\right) =\vec{x}_{0}^{\mathrm{T}}P_{R}\vec{x}_{0}=%
\breve{x}_{0}^{\mathrm{H}}P_{C}\breve{x}_{0}.  \label{min}
\end{equation}
\end{remark}

\begin{remark}
If Theorem \ref{th3} is applied on the normal linear system (\ref{normal}), $%
Q_{2}=0,$ and $R_{2}=0$, it will be reduced to its standard LQR solution. In
fact, as all the known matrices in (\ref{are1}) are diagonal, a trivial
diagonal solution $P_{C}=\mathrm{diag}\{P_{0},P_{0}^{\#}\}>0$ exists
(namely, $\left \{ P_{1},P_{2}\right \} =\left \{ P_{0},0\right \} $), where
$P_{0}$ solves the normal ARE%
\begin{equation}
-Q_{1}=\left \{
\begin{array}{ll}
A_{1}^{\mathrm{H}}P_{0}+P_{0}A_{1}-P_{0}B_{1}R^{-1}B_{1}^{\mathrm{H}}P_{0},
& t\in \mathbf{R}^{+}, \\
A_{1}^{\mathrm{H}}P_{0}A_{1}-P_{0}-A_{1}^{\mathrm{H}%
}P_{0}B_{1}S_{0}^{-1}B_{1}^{\mathrm{H}}P_{0}A_{1}, & t\in \mathbf{Z}^{+},%
\end{array}%
\right.  \label{eq85}
\end{equation}%
in which $S_{0}=R_{1}+B_{1}^{\mathrm{H}}P_{0}B_{1}.$ On the other hand, as
the system is stabilizable and $Q_{1}>0$, the positive definite solution to (%
\ref{eq85}) is unique \cite{kucera73kyb}, \cite{kucera72kyb}, %
\cite{lr95book}. Thus $\left \{ P_{1},P_{2}\right \} =\left \{
P_{0},0\right \} $ is the unique solution to (\ref{eqareA})--(\ref{eqareB})
and it follows from (\ref{eqgain}) that $K_{2}^{\ast }=0,$ namely, the
optimal controller is exactly the normal linear feedback (\ref%
{normalfeedback}) with $K_{1}=K_{1}^{\ast }$ given by (\ref{eqgain}). This
coincides with the well-known LQR theory \cite{lr95book}.
\end{remark}

If we apply Theorem \ref{th3} on the antilinear system (\ref{antilinear}),
we obtain the following corollary.

\begin{corollary}
\label{coro4}Consider the antilinear system (\ref{antilinear}) which is
assumed to be stabilizable, and $R_{2}=0$.

\begin{enumerate}
\item When $t\in \mathbf{R}^{+}$, there is a unique bimatrix $\left \{
P_{1},P_{2}\right \} >0$ to
\begin{equation}
\left \{
\begin{array}{rll}
-Q_{1} & = & A_{2}^{\mathrm{H}}P_{2}+P_{2}^{\#}A_{2}-P_{2}^{%
\#}B_{2}R_{1}^{-1}B_{2}^{\mathrm{H}}P_{2}-P_{1}B_{2}^{\#}R_{1}^{-\#}B_{2}^{%
\mathrm{T}}P_{1}, \\
-Q_{2} & = & A_{2}^{\mathrm{T}}P_{1}+P_{1}^{\#}A_{2}-P_{1}^{%
\#}B_{2}R_{1}^{-1}B_{2}^{\mathrm{H}}P_{2}-P_{2}B_{2}^{\#}R_{1}^{-\#}B_{2}^{%
\mathrm{T}}P_{1}.%
\end{array}%
\right.  \label{are7}
\end{equation}%
Moreover, the optimal control is the full state feedback (\ref{optcontrol})
with%
\begin{equation}
K_{1}^{\ast }=-R_{1}^{-1}B_{2}^{\mathrm{H}}P_{2},\;K_{2}^{\ast
}=-R_{1}^{-\#}B_{2}^{\mathrm{T}}P_{1}.  \label{eqk1k21}
\end{equation}

\item When $t\in \mathbf{Z}^{+}$ and $Q_{2}=0$, if
\begin{equation}
-Q_{1}=A_{2}^{\mathrm{H}}P_{1}^{\#}A_{2}-A_{2}^{\mathrm{H}%
}P_{1}^{\#}B_{2}S_{1}^{-1}B_{2}^{\mathrm{H}}P_{1}^{\#}A_{2}-P_{1},
\label{are2}
\end{equation}%
where $S_{1}=R_{1}+B_{2}^{\mathrm{H}}P_{1}^{\#}B_{2},$ has a solution $%
P_{1}>0$, then such a solution is unique, and the optimal controller is the
normal state feedback (\ref{normalfeedback}) with $K_{1}=K_{1}^{\ast }$
defined by%
\begin{equation}
K_{1}^{\ast }=-\left( R_{1}+B_{2}^{\mathrm{H}}P_{1}^{\#}B_{2}\right)
^{-1}B_{2}^{\mathrm{H}}P_{1}^{\#}A_{2},  \label{eqk0}
\end{equation}%
and the optimal value of $J\left( u\right) $ given by $J_{\min }\left(
u\right) =x_{0}^{\mathrm{H}}P_{1}x_{0}.$
\end{enumerate}

Moreover, in both cases, the closed-loop system is asymptotically stable.
\end{corollary}

The point of Corollary \ref{coro4} is that, for system (\ref{antilinear}),
full state feedback is necessary with $t\in \mathbf{R}^{+}$, while the
normal state feedback is sufficient when $t\in \mathbf{Z}^{+}$. Item 2 also
improves a little the main results in \cite{wqls16jfi} where the system is
assumed to be controllable. Another interesting point, though not proven, is
that the so-called anti-ARE (\ref{are2}) may always have a positive definite
solution if the system is stabilizable. Numerical examples support this
conjecture, yet a proof is not available at present.

Solutions to the non-standard anti-ARE (\ref{are2}) have been investigated
in our early papers \cite{lzl14amc} and \cite{zcl13amc}.

\subsection{State observer design}

We finally discuss briefly in this subsection the state observer design for
system (\ref{sys}). For simplicity, we let $\left \{ D_{1},D_{2}\right \} =%
\mathcal{O}_{p\times m}$. Motivated by the state feedback design, we present
the following full-order state observer%
\begin{equation}
z^{+}=\left \{ A_{1},A_{2}\right \} z+\left \{ B_{1},B_{2}\right \} u+\left
\{ L_{1},L_{2}\right \} \left( \left \{ C_{1},C_{2}\right \} z-y\right) ,
\label{observer}
\end{equation}%
where $\left \{ L_{1},L_{2}\right \} \in \{ \mathbf{C}^{n\times p},\mathbf{C}%
^{n\times p}\}$ is an observer gain bimatrix to be designed, and $z\in
\mathbf{C}^{n}$ is the observer state. Then, by denoting $e=x-z,$ we obtain
the observer-error system%
\begin{equation}
e^{+}=\left( \left \{ A_{1},A_{2}\right \} +\left \{ L_{1},L_{2}\right \}
\left \{ C_{1},C_{2}\right \} \right) e.  \label{error}
\end{equation}%
Taking $\vec{\cdot}$ on both sides of (\ref{error}) gives%
\begin{equation*}
\vec{e}^{+}=\left( \left \{ A_{1},A_{2}\right \} _{\circ }+\left \{
L_{1},L_{2}\right \} _{\circ }\left \{ C_{1},C_{2}\right \} _{\circ }\right)
\vec{e}.
\end{equation*}%
Then, if there exists a real gain $L\in \mathbf{R}^{2n\times 2p},$ such that
$\left \{ A_{1},A_{2}\right \} _{\circ }+L\left \{ C_{1},C_{2}\right \}
_{\circ } $ is asymptotically stable (possesses a prescribed eigenvalue set $%
\Gamma $ that is symmetric with respect to the real axis), then $%
\lim_{t\rightarrow \infty }\left \Vert \vec{e}(t)\right \Vert
=\lim_{t\rightarrow \infty }\left \Vert x(t)-z(t)\right \Vert =0,$ which
means that $z$ converges asymptotically to the true state of system, namely,
(\ref{observer}) is a state observer. Moreover, $\left \{
L_{1},L_{2}\right
\} $ can be determined by%
\begin{equation}
\left[
\begin{array}{c}
L_{1} \\
L_{2}%
\end{array}%
\right] =H_{n}LH_{p}^{\mathrm{H}}\left[
\begin{array}{c}
I_{p} \\
0_{p\times p}%
\end{array}%
\right] .  \label{eql1l2}
\end{equation}%
Finally, there exists a real gain $L\in \mathbf{R}^{2n\times 2p},$ such that
$\left \{ A_{1},A_{2}\right \} _{\circ }+L\left \{ C_{1},C_{2}\right \}
_{\circ } $ is asymptotically stable (possesses a prescribed eigenvalue set $%
\Gamma $) if and only if $(\left \{ A_{1},A_{2}\right \} _{\circ },\left \{
C_{1},C_{2}\right \} _{\circ })$ is detectable (observable), which is
further equivalent to the detectability (observability) of $(\left \{
A_{1},A_{2}\right \} _{\diamond },\left \{ C_{1},C_{2}\right \} _{\diamond
}),$ namely,%
\begin{equation*}
\mathrm{rank}\left[
\begin{array}{cc}
sI_{n}-A_{1} & -A_{2}^{\#} \\
-A_{2} & sI_{n}-A_{1}^{\#} \\
C_{1} & C_{2}^{\#} \\
C_{2} & C_{1}^{\#}%
\end{array}%
\right] =2n,\; \forall s\in \mathbf{D}.
\end{equation*}%
We call that system (\ref{sys}) is detectable if $(\left \{
A_{1},A_{2}\right \} _{\circ },\left \{ C_{1},C_{2}\right \} _{\circ })$ is
detectable.

We can also state a separation principle in the observer-based state
feedback scheme, say, the observer gain bimatrix $\left \{ L_{1},L_{2}\right
\} $ and the state feedback gain bimatrix $\left \{ K_{1},K_{2}\right \} $
can be designed separately. Finally, we point out that the reduced-order and
functional observers design can be considered in a similar way and are
omitted for brevity.

\section{\label{sec5}Concluding remarks}

This paper has investigated a class of $n$ dimensional complex-valued linear
systems. By using the the so-called bimatrix and its properties, and
representing the considered linear system as a $2n$ dimensional real-valued
linear system and a $2n$ dimensional complex-valued linear system, the
solutions, controllability, observability and stability are studied and
explicit solutions/conditions were provided in terms of bimatrices of the
original system. Stabilization, eigenvalue assignment, linear quadratic
regulation (LQR), and observer design problems for the complex-valued linear
system were solved and conditions were also expressed in terms of the
original system parameters. Particularly, the stability analysis and
solutions to LQR problem are characterized by the so-called Lyapunov
bimatrix equation and the bimatrix algebraic Riccati equation, which are
equivalent to coupled matrix equations. The obtained results are degenerated
to the existing ones when the system is reduced to the normal linear system,
and improve the existing ones when the system is reduced to the antilinear
system.

We may consider more general models than (\ref{sys}). One case is allowing
time-varying parameters, and the other is allowing additional terms on the
left hand side of (\ref{sys}), namely,%
\begin{equation}
\left \{
\begin{array}{rl}
\left \{ E_{1},E_{2}\right \} x^{+} & =\left \{ A_{1},A_{2}\right \} x+\left
\{ B_{1},B_{2}\right \} u, \\
y & =\left \{ C_{1},C_{2}\right \} x+\left \{ D_{1},D_{2}\right \} u.%
\end{array}%
\right.  \label{sys4}
\end{equation}%
If $\{E_{1},E_{2}\}$ is nonsingular, this is equivalent to (\ref{sys2})
where $\{A_{1},A_{2}\}$ and $\{B_{1},B_{2}\}$ are replaced by $%
\{E_{1},E_{2}\}^{-1}\{A_{1},A_{2}\}$ and $\{E_{1},E_{2}\}^{-1}\{B_{1},B_{2}%
\},$ respectively. If $\{E_{1},E_{2}\}$ is not nonsingular, system (\ref%
{sys4}) should be studied in the descriptor system framework %
\cite{duan10book}. The approaches and results in this paper can also be
extended to other more general complex-valued systems, for example,
complex-valued nonlinear systems, complex-valued linear systems in matrix
forms (see \cite{zdz09ijss} for normal matrix-valued linear systems), and
complex-valued time-delay systems.

\end{document}